\documentclass{article}

\usepackage{amssymb}
\usepackage{amsmath}
\usepackage{amsthm}
\usepackage{mathtools}
\usepackage{graphicx}
\usepackage{epstopdf}
\usepackage{caption}
\usepackage{multirow}
\usepackage{booktabs}
\usepackage[flushleft]{threeparttable}
\usepackage{adjustbox}
\usepackage{array}
\usepackage{bm}
\usepackage[svgnames]{xcolor}
\usepackage{url}
\usepackage[colorlinks=true, linkcolor=Brown, citecolor=DarkCyan, urlcolor=Olive]{hyperref}
\usepackage{array}
\newcolumntype{C}{>{$}c<{$}}
\usepackage{appendix}
\usepackage[square, numbers]{natbib}

\usepackage{xcolor,colortbl}
\definecolor{Gray}{gray}{0.85}
\definecolor{LightGray}{gray}{0.95}
\newcolumntype{a}{>{\columncolor{Gray}}c}
\newcolumntype{b}{>{\columncolor{LightGray}}c}

\newcommand{\vecn}{{n}}
\newcommand{\vecu}{{u}}
\newcommand{\vecv}{{v}}
\newcommand{\vecw}{{w}}
\newcommand{\vecx}{{x}}

\newcommand{\gothg}{\mathfrak{g}}
\newcommand{\gothf}{\mathfrak{f}}
\newcommand{\gothp}{\mathfrak{p}}
\newcommand{\gothr}{\mathfrak{r}}

\newcommand{\gothw}{\mathfrak{w}}
\newcommand{\gothu}{\mathfrak{u}}
\newcommand{\gothv}{\mathfrak{v}}

\newcommand{\gothq}{\mathfrak{q}}

\newcommand{\mixperm}{\mathfrak{K}\,}
\newcommand{\gothpp}{\tilde{\mathfrak{p}}_h}

\newcommand{\calC}{\mathcal{C}}

\newcommand{\Tau}{\mathcal{T}}
\newcommand{\calK}{{\mathcal{K}}}

\newcommand{\trace}{\mathrm{tr}\,}
\newcommand{\mixdiv}{\mathfrak{D}\cdot}
\newcommand{\mixgrad}{\mathbb{D}\,}
\newcommand{\mixTraceN}[1]{\mathfrak{T}_N{#1}}
\newcommand{\mixTraceD}[1]{\mathfrak{T}_D{#1}}
\newcommand{\mixTraceX}[1]{\mathfrak{T}_X{#1}}

\newcommand{\mixLIIscalar}{L^2\left(\Omega\right)}

\newcommand{\mixHI}{H^1(\Omega)}
\newcommand{\mixHIzero}{H^1_0(\Omega)}

\newcommand{\mixHdiv}{H(\mathrm{div};\Omega,\Gamma)}
\newcommand{\mixHdivU}{H(\mathrm{div};\Omega,\Gamma;U)}
\newcommand{\mixHdivUa}[1]{H(\mathrm{div};\Omega,\Gamma;U_{#1})}
\newcommand{\mixHdivX}[1]{H(\mathrm{div};\Omega,\Gamma;{#1})}

\newcommand{\shigh}{\hat{\mathcal{S}}_i}
\newcommand{\slow}{\check{\mathcal{S}}_i}

\newcommand{\mhalf}{-{\frac{1}{2}}}
\newcommand{\phalf}{{\frac{1}{2}}}

\newcommand{\hatj}{{\hat{\jmath}}}
\newcommand{\checkj}{{\check{\jmath}}}

\newcommand{\ssum}{{\textstyle{\sum}}}

\newcommand{\inner}[2]{{\left\langle{#1},{#2}\right\rangle}}
\newcommand{\innerangle}[2]{{\left\langle{#1},{#2}\right\rangle}}

\DeclarePairedDelimiter\abs{\lvert}{\rvert}%
\DeclarePairedDelimiter\norm{\lVert}{\rVert}%
\makeatletter
\let\oldabs\abs
\def\abs{\@ifstar{\oldabs}{\oldabs*}}
\let\oldnorm\norm
\def\norm{\@ifstar{\oldnorm}{\oldnorm*}}
\makeatother

\newcommand{\tnorm}[1]{{\left\vert\kern-0.25ex\left\vert\kern-0.25ex\left\vert #1 
    \right\vert\kern-0.25ex\right\vert\kern-0.25ex\right\vert}}
\newcommand{\tnormstar}[1]{\tnorm{#1}_*}

\newcommand{\fullnorm}[1]{{\left\vert\kern-0.25ex\left\vert\kern-0.25ex\left[ #1 
    \right]\kern-0.25ex\right\vert\kern-0.25ex\right\vert}}

\newtheorem{theorem}{Theorem}

\newtheorem{lemma}{Lemma}

\newtheorem{definition}{Definition}

\newtheorem{remark}{Remark}

\begin{document}
\title{\textit{A posteriori} error estimates for hierarchical mixed-dimensional elliptic equations}

\author{Jhabriel Varela\footnote{Corresponding author e-mail: \href{mailto:jhabriel.varela@uib.no}{jhabriel.varela@uib.no.}} $^{,}$\footnote{Center for Modeling of Coupled Subsurface Dynamics, University of Bergen, P.O. Box 7800, N-5020 Bergen, Norway.} \and Elyes Ahmed\footnote{SINTEF Digital,  Mathematics  and  Cybernetics,  P.O.  Box 124 Blindern, N-0314 Oslo, Norway.} \and Eirik Keilegavlen$^{\dagger}$ \and Jan Martin Nordbotten$^{\dagger}$ \and Florin Adrian Radu$^{\dagger}$}

\date{\today}

\maketitle

\begin{abstract}
Mixed-dimensional elliptic equations exhibiting a hierarchical structure are commonly used to model problems with high aspect ratio inclusions, such as flow in fractured porous media. We derive general abstract estimates based on the theory of functional \textit{a posteriori} error estimates, for which guaranteed upper bounds for the primal and dual variables and two-sided bounds for the primal-dual pair are obtained. We improve on the abstract results obtained with the functional approach by proposing four different ways of estimating the residual errors based on the extent the approximate solution has conservation properties, i.e.: (1) no conservation, (2) subdomain conservation, (3) grid-level conservation, and (4) exact conservation. This treatment results in sharper and fully computable estimates when mass is conserved either at the grid level or exactly, with a comparable structure to those obtained from grid-based \textit{a posteriori} techniques. We demonstrate the practical effectiveness of our theoretical results through numerical experiments using four different discretization methods for synthetic problems and applications based on benchmarks of flow in fractured porous media.\vspace{3mm} \\
\textbf{Keywords}: mixed-dimensional geometry, functional a posteriori error estimates, fractured porous media \vspace{3mm}\\
\textbf{Classification}: 65N15, 76S05, 35Q86 
\end{abstract}

\newpage
\section{Introduction}

Mixed-dimensional partial differential equations (mD-PDEs) arise when partial differential equations interact on domains of different topological dimensions~\cite{nordbotten2019mixed}. Prototypical examples include models of thin inclusions in elastic materials~\cite{antman1995nonlinear,ciarlet1997mathematical,boon2020stable}, blood flow in human vasculature \cite{dangelo2008coupling,koppl2015local,hodneland2019new}, root water uptake systems~\cite{koch2018new}, and flow in fractured porous media~\cite{alboin2002modeling,formaggia2014reduced,ahmed2017reduced}. The latter example has an appealing mathematical structure, in that the model equations allow for a hierarchical representation where each subdomain (matrix, fractures, fracture intersections, and intersection points) only has direct interaction with subdomains of topological dimension one higher or one lower~\cite{boon2018robust}. Such hierarchical mD-PDEs are the topic of the current paper. 

\begin{figure}
    \includegraphics[width=.99\textwidth]{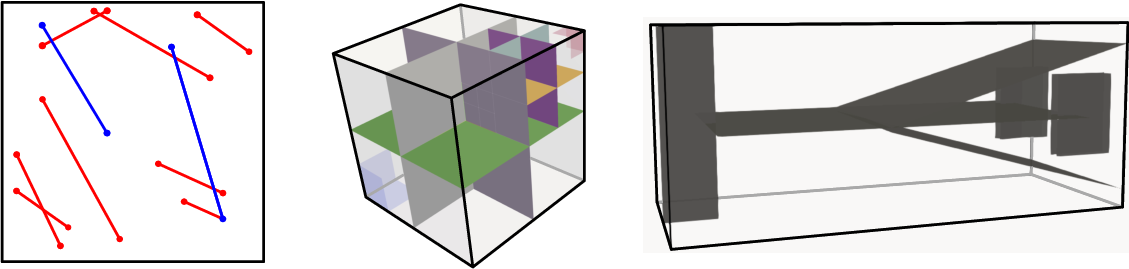}
    \caption{Example geometries falling within the context of hierarchical mixed-dimensional geometries studied herein. Left figure corresponds to a 2d benchmark problem \cite{flemisch2018benchmarks} while the two remaining correspond to 3d benchmark problems \cite{flemisch2020verification}.}
    \label{fig:md-examples}
\end{figure}

mD-PDEs are intrinsically linked to the underlying geometric representation, which, in a certain sense, generalizes the usual notion of the domain. One can then define sets of suitable functions (and function spaces) on this geometry, and these sets are then naturally interpreted as mixed-dimensional (mD) functions. Exploiting this concept, one can generalize the standard differential operators to mappings between mD functions and thus obtain an mD calculus. The fact that this mD calculus inherits standard properties of calculus, particularly partial integration (relative to suitable inner products), a de~Rham complex structure, and a Poincaré-Friedrichs inequality, was recently established using the language of exterior calculus on differential forms~\cite{boon2020functional}. 

The inherent geometric generality of hierarchical mD-PDEs also demand the same level of abstraction of \textit{a posteriori} error estimation techniques. This requirement makes error estimates of the functional type particularly well-suited for the task~\cite{repin2000variational, repin2003twosided, neittaanm2004reliable, repin2007mixed, repin2008posteriori, pauly2020solution}. The most attractive feature of this approach is that error estimates are derived using purely functional methods~\cite{repin2008posteriori}. The bounds are therefore agnostic to the way approximated solutions are obtained in the energy space, and the only undetermined constants arise from Poincar{\'e}-type inequalities~\cite{kurzFunctional2021}. 

However, unlike other types of error estimates~\cite{verfurth1999residual,zz1987estimator,oden2001Goal, vohralik2010unified,ainsworth2007mixed}, this generality makes standard functional estimates of limited applicability to hierarchical elliptic mD-PDEs due to the following reasons: (1) for general fracture networks, the mixed-dimensional Poincar{\'e} constant is not easily computable, and (2) since Poincar{\'e} constants are proportional to the diameter of the physical domain, residual estimators cannot exhibit superconvergent properties.

To circumvent the aforementioned issues, we exploit the fact that Poincar{\'e}-type inequalities imply weighted norms~\cite{pechstein2013weighted, rathmair2019poincare}, and use spatially-dependent weights to control the residual norms. We show both theoretically and numerically that this treatment leads to sharper estimates when approximations to the exact solution satisfy mass conservation in a given partition of the domain.

In view of the preceding discussion, our aim is therefore to obtain \textit{a posteriori} error estimates for the approximate solution to the mD~scalar elliptic equation~\cite{boon2018robust, nordbotten2019unified, boon2020functional}, where the mD~Laplace equation for geometries such as those illustrated in Figure~\ref{fig:md-examples} is described in detail in Section~\ref{sec:extensionToFracNet}.

We remark that while a broad range of \textit{a~posteriori} error techniques are available for mono-dimensional problems, existing error bounds for mD~models are far more scarce. Moreover, the ones available, are restricted to specific cases (e.g., in the context of mortar methods \cite{wohlmuth1999mortar, belhachmi2003mortar, wheeler2005mortar, pencheva2013mortar} and fractured porous media \cite{chen2017fractured,mghazli2019fractured,hecht2019residual}) with far less geometric generality than what we present here. Thus, for practical problems, \textit{a posteriori} error bounds for mD~geometries have until now essentially not been available.

The rest of the paper is structured as follows: Section~\ref{sec:singleFrac} is devoted to introducing the model problem, functional spaces, and variational formulations for the case of a single 1d fracture embedded in a 2d matrix. The section is concluded by providing a first upper bound for the primal variable. In Section~\ref{sec:extensionToFracNet}, we generalize the results from Section~\ref{sec:singleFrac} to the case of fracture networks and introduce the necessary tools to perform the \textit{a posteriori} analysis in an mD~setting. After reviewing necessary tools from functional analysis in Section~\ref{sec:functools}, in Section~\ref{sec:posteriori}, we provide our main results starting from a generic abstract estimate and then considering specific cases depending upon the degree of accuracy at which residual terms are approximated. In Section~\ref{sec:concreteBounds}, we introduce the approximated problem using mixed-finite element methods and thus make the estimates concrete. Sections~\ref{sec:num_val}~and~\ref{sec:numex_app} deal, respectively, with numerical validations and practical applications of the derived bounds. Finally, in Section~\ref{sec:conclusion}, we present our concluding remarks.

\section{Upper bounds for a single fracture \label{sec:singleFrac}}

In this section, we introduce the model problem together with functional spaces and the variational formulations for the case of a single 1d line embedded in a 2d matrix, as illustrated in Figure~\ref{fig:singleFracture}. Furthermore, a first upper bound for the primal variable is derived following the classical functional approach. We remark that the case of a single fracture embedded in a matrix has been analyzed before. For example, \cite{chen2017fractured} and \cite{hecht2019residual} proposed error estimators based on the residual approach, whereas \cite{mghazli2019fractured} obtained guaranteed a posteriori error estimates using the approach of Vohralík \cite{vohralik2010unified}.

\subsection{The model problem for a single fracture}
\label{ref:mod_single}

\begin{figure}[tbp]
    \centerline{\includegraphics[height=0.225\textheight]{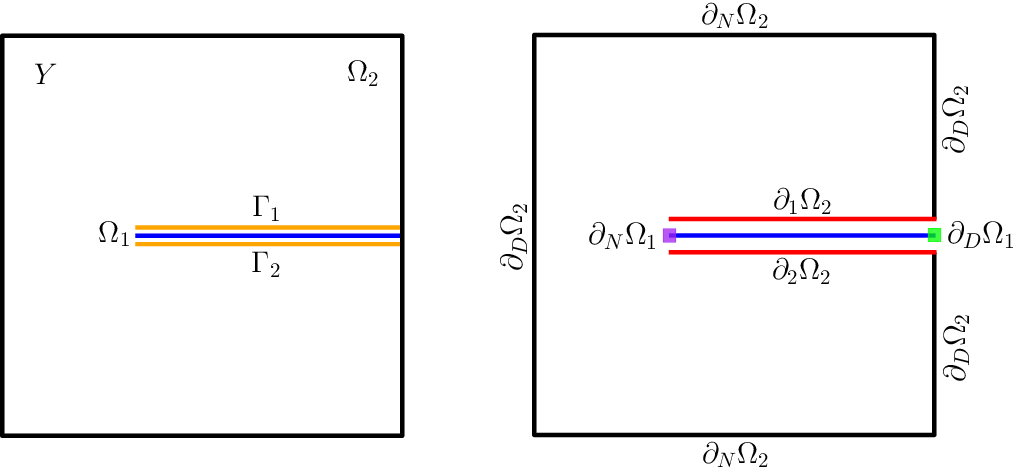}}
	\caption{A horizontal 1d fracture embedded in a 2d matrix. Left: Subdomains and interfaces. Right: Boundary conditions. For the fracture, the purple square denotes a no-flux boundary condition, whereas the green square a Dirichlet boundary condition. Note that $\partial_1\Omega_2$, $\Gamma_1$, $\Omega_1$, $\Gamma_2$, $\partial_2\Omega_2$, all coincide spatially. For illustrative purposes, however, they are placed in different locations.  \label{fig:singleFracture}}
\end{figure}

Before writing the set of equations describing general fracture networks, let us first introduce the governing equations of a simpler configuration; that is, a unit square domain $Y\subset \mathbb{R}^2$ decomposed as a 1d fracture $\Omega_1$ embedded in a 2d matrix $\Omega_2$ as shown in the left Figure~\ref{fig:singleFracture}. Interfaces $\Gamma_1$ and $\Gamma_2$, at each side of $\Omega_1$, establish the link between $\Omega_2$ and $\Omega_1$. The model presented below is well-established for these problems, and we point the reader to the references for further justification of this system~\cite{martin2005modeling, boon2018robust, nordbotten2019unified}. 

The strong form of the governing equations in $\Omega_2$ reads
\begin{subequations}
    \begin{alignat}{2}
        \nabla \cdot \vecu_2 &= f_2, \qquad&&\mathrm{in}\,\Omega_2,  \label{eq:massConservationMatrix} \\
        \vecu_2 &= - \mathcal{K}_2 \nabla\, p_2, \qquad&&\mathrm{in}\,\Omega_2, \label{eq:DarcyLawMatrix} \\
        \vecu_2\cdot \vecn_2 &= \lambda_1, \qquad&&\mathrm{on}\,\partial_1\Omega_2, \label{eq:InternalBoundary1} \\
        \vecu_2\cdot \vecn_2 &= \lambda_2, \qquad&&\mathrm{on}\,\partial_2\Omega_2, \label{eq:InternalBoundary2} \\
        \vecu_2\cdot\vecn_2 &= 0, \qquad&&\mathrm{on}\,\partial_N\Omega_2,\label{eq:NeumannBoundaryMatrix} \\
        p_2 &= g_{D,2}, \qquad&&\mathrm{on}\,\partial_D\Omega_2. \label{eq:DirichletBoundaryMatrix}
    \end{alignat}
\end{subequations}
Here, \eqref{eq:massConservationMatrix} is the mass conservation equation, $\vecu_2$ is the matrix velocity, and $f_2$ an external source. The fluid velocity is given by the standard Darcy's law~\eqref{eq:DarcyLawMatrix}, where $\mathcal{K}_2$ is the matrix permeability; a bounded, symmetric, and positive-definite $2\times 2$ tensor, and $p_2$ is the fluid pressure. 

Equations \eqref{eq:InternalBoundary1} and \eqref{eq:InternalBoundary2}  require that at each side of the internal boundary of $\Omega_2$, the normal component of $\vecu_2$ to match the interface (mortar) fluxes $\lambda_1$ and $\lambda_2$. To fix the direction of the normal vector on internal boundaries, we require $\vecn_2$ pointing from the higher- to the lower-dimensional subdomain. No flux conditions are prescribed in \eqref{eq:NeumannBoundaryMatrix}, where $\vecu_2\cdot\vecn_2$ represents the outer normal flux across $\partial_N\Omega_2$. Finally, Dirichlet boundary conditions are imposed in \eqref{eq:DirichletBoundaryMatrix}, where $g_{D,2}$ is a prescribed function on the Dirichlet boundary. 

In the fracture $\Omega_1$, the equations are given by
\begin{subequations}
\begin{alignat}{2}
    \nabla_{1} \cdot\vecu_1 - \left( \lambda_1 +  \lambda_2\right) &= f_1, \qquad&&\mathrm{in}\,\Omega_1,\label{eq:massConservationFracture}\\
    \vecu_1 &= - \mathcal{K}_1 \nabla_{1}\, p_1, \qquad&&\mathrm{in}\,\Omega_1,\label{eq:DarcyLawFracture} \\
    \vecu_1 \cdot \vecn_1 &= 0, \qquad&&\mathrm{on}\,\partial_N\Omega_1,\label{eq:NeumannBoundaryFracture} \\
    p_1 &= g_{D,1}, \qquad&&\mathrm{on}\,\partial_D\Omega_1.\label{eq:DirichletBoundaryFracture}
\end{alignat}
\end{subequations}
In~\eqref{eq:massConservationFracture}, $\nabla_1\cdot(\cdot)=\frac{d}{dx}(\cdot)=\nabla_1(\cdot)$ are the divergence and gradient operators acting in the tangent space of $\Omega_1$, $\vecu_1$ is the tangential fracture velocity, the term in parentheses represents the jump in normal fluxes from the adjacent interfaces $\Gamma_1$ and $\Gamma_2$ onto $\Omega_1$, and $f_1$ is an external source.  

The tangential velocity $u_1$ is again expressed via Darcy's law~\eqref{eq:DarcyLawFracture}, where in a slight abuse of notation, we use $\mathcal{K}_1$ to refer to the tangential component of the fracture permeability, which is again assumed to be positive and bounded from above. Finally, \eqref{eq:NeumannBoundaryFracture} and \eqref{eq:DirichletBoundaryFracture} are the Neumann and Dirichlet boundary conditions, respectively. Again, we use $g_{D,1}$ to denote a prescribed function on the Dirichlet part of the fracture boundary.

To close the system of equations, we must specify a constitutive relationship for the interface fluxes. Here, we use a Darcy-type law~\cite{martin2005modeling}, where mortar fluxes are linearly related to pressure jumps
\begin{subequations}
\begin{alignat}{2}
    \lambda_1 &= -\kappa_1 \left(p_1 - p_2\right), \qquad&&\mathrm{on}\,\Gamma_1,\label{eq:DarcyLawInterface1} \\
    \lambda_2 &= -\kappa_2 \left(p_1 - p_2\right), \qquad&&\mathrm{on}\,\Gamma_2,\label{eq:DarcyLawInterface2}
\end{alignat}    
\end{subequations}
with $\kappa_1$ and $\kappa_2$ representing the effective normal permeability on $\Gamma_1$ and $\Gamma_2$, respectively. We restrict our analysis to the case where $\kappa_1$ and $\kappa_2$ are non-degenerate. Thus, following~\cite{boon2018robust}, we further require the existence of two constants $\gamma_1$ and $\gamma_2$ such that $0 < \gamma_1 \leq \kappa^{-1}_j \leq \gamma_2 < \infty$ for $j\in\{1,2\}$.

\subsection{Functional spaces and variational formulations}
\label{sec:func_simple}
Let us now present the primal weak formulation of the single fracture model from the previous section. To this aim, consider first the energy space with vanishing traces on Dirichlet boundaries 
\begin{align}
    H_{0}^1(\Omega_i) &= \{q_i\in H^1(\Omega_i) : \mathrm{tr}_{\partial_D\Omega_i}~q_i = 0 \},
\end{align}
and the product spaces
\begin{align}
   H^1(\Omega) = H^1(\Omega_1)\times H^1(\Omega_2) \qquad\mathrm{and}\qquad H_{0}^1(\Omega) = H_{0}^1(\Omega_1)\times H_{0}^1(\Omega_2).
   \label{eq:prod2}
\end{align}

Furthermore, let $\inner{\cdot}{\cdot}_{\Omega_i}$ and $\innerangle{\cdot}{\cdot}_{\Gamma_j}$ denote respectively the $L^2$--inner products on $\Omega_i$ and $\Gamma_j$, and $\norm{\cdot}_{\Omega_i}$ and $\norm{\cdot}_{\Gamma_j}$ the relevant $L^2$--norms. Finally, we denote by $g=[g_{1}, g_{2}]\in H^1(\Omega)$ two functions extending the boundary data into the domains, and thus satisfying $\mathrm{tr}_{\partial_D\Omega_i} g_i = g_{D,i}$. We now state the primal weak problem as:

\begin{definition}[Primal weak formulation for a single fracture\label{def:primalWeakDefinition}] Let $p=[p_1,p_2]$ and $g=[g_{1}, g_{2}]\in H^1(\Omega)$. Then find $p\in H^1_0(\Omega)+g$ such that 
\begin{align}
    \sum_{i=1}^2 \inner{\mathcal{K}_i\,\nabla_i\,p_i}{\nabla_i\,q_i}_{\Omega_i} + \sum_{j=1}^2 \innerangle{\kappa_j\left(p_1 - \mathrm{tr}_{\partial_j\Omega_2}\,p_2 \right)}{q_1 - \mathrm{tr}_{\partial_j\Omega_2}\,q_2}_{\Gamma_j} \nonumber \\
    = \sum_{i=1}^2 \inner{f_i}{q_i}_{\Omega_i}, \quad\forall\,q = [q_1, q_2] \in H_{0}^1(\Omega).\label{eq:primalSingleFracture}
\end{align}
\end{definition}

Refer to Appendix~\ref{sec:primalWeakFormProof} for the derivation of the primal weak form from the strong form in Section~\ref{ref:mod_single}. We see directly from equation \eqref{eq:primalSingleFracture} that the primal weak form has a minimization structure subject to the stated conditions on $\mathcal{K}_i$ and $\kappa_j$, and well-posedness follows by standard arguments. 

A dual weak form for the model problem, with explicit representation of the subdomain velocities and mortar fluxes, can also be constructed. We first define the space $H(\mathrm{div}; \Omega_i, \partial_X \Omega)$ as the space of $L^2$-vector functions on $\Omega_i$ with weak divergence in $L^2(\Omega_i)$ and zero trace on the part of the boundary indicated by $\partial_X \Omega$. Then, we denote the product spaces of $H(\mathrm{div})$-functions that are zero on Neumann, and on Neumann and internal boundaries as:
\begin{align}
{V} &= H(\mathrm{div}; \Omega_1, \partial_N \Omega_1)\times H(\mathrm{div}; \Omega_2, \partial_N \Omega_2), \\
{V}_0 &= H(\mathrm{div};\Omega_1,\partial_N\Omega_1) \times  H(\mathrm{div};\Omega_2,\partial_N \Omega_2\cup\partial_1\Omega_2\cup\partial_2\Omega_2).
\end{align}

Furthermore, we define the $L^2$-product spaces on the domains:
\begin{equation} 
L^2(\Omega) =  L^2(\Omega_1)\times L^2(\Omega_2), \qquad
L^2(\Gamma) =  L^2(\Gamma_1)\times L^2(\Gamma_2).
\end{equation}
With these spaces in hand, we consider the standard linear extension operators from internal boundaries onto domains denoted $\mathcal{R}_j := L^2(\Gamma_j) \to H(\mathrm{div}; \Omega_2, \partial_N \Omega_2)$, such that $\mathcal{R}_j$ satisfies for all $\lambda_j \in L^2(\Gamma_j)$
\begin{equation}
\label{eq:reconstruction}
    \mathrm{tr}_{\partial_j\Omega_2}\,(\mathcal{R}_j\,\lambda_j) \cdot \vecn_2 = 
    \begin{cases}
         \lambda_j & \mathrm{on}\,\partial_j\Omega_2 \\
        0 & \mathrm{on}\,\partial\Omega\setminus \partial_j\Omega_2
    \end{cases}.
\end{equation}
The precise choice of the extension operator $\mathcal{R}_j$ is not important; however, the natural choice based on the solution of an auxiliary elliptic equation is reasonable~\cite{boon2018robust}. We naturally extend the definition of $\mathcal{R}_j$ to $\mathcal{R} := L^2(\Gamma) \to V$ by requiring that for $[\lambda_1,\lambda_2] \in L^2(\Gamma)$, then $[\vecu_1, \vecu_2]=\mathcal{R}\lambda$ satisfies $\vecu_1=0$ and $\vecu_2=\mathcal{R}_1 \lambda_1 + \mathcal{R}_2 \lambda_2$.

The above constructions allow us to represent subdomain fluxes as
\begin{equation}
 \vecu = \vecu_0 +  \mathcal{R}\lambda,
 \label{eq:decompose_flux}
\end{equation}
where $\vecu_0 \in {V}_0$ and  $\lambda \in L^2(\Gamma)$. This motivates the construction of a compound $H(\mathrm{div})$-type spaces, as
\begin{equation}
    \mixHdiv = {V}_{0} \times L^2(\Gamma).\label{eq:firstHdiv}
\end{equation}
This construction will become key when we generalize to more complex geometries in the next section.

\begin{remark}[On the regularity of $\mixHdiv$] It is worth remarking that the restriction of space $\mixHdiv$ to the domain $\Omega_2$ has slightly enhanced regularity relative to the standard space ${H}(\mathrm{div};\Omega_2)$, as this latter space has normal traces which do not lie in $L^2(\Gamma_1)$ nor $L^2(\Gamma_2)$.
\end{remark}

\begin{definition}[Dual weak formulation for a single fracture.\label{def:dualSingleFrac}] Let $\vecu_0 = [\vecu_{0,1},\vecu_{0,2}]$, $\lambda = [\lambda_1,\lambda_2]$, $p = [p_1, p_2]$. Then find $(\vecu_0,\lambda,p) \in \mixHdiv \times L^2(\Omega)$ such that 
\begin{subequations}
\begin{alignat}{2}
    &\inner{\mathcal{K}_2^{-1}\left(\vecu_{0,2} + \mathcal{R}_1\lambda_1 + \mathcal{R}_2\lambda_2\right)}{\vecv_{0,2}}_{\Omega_2} + \inner{\mathcal{K}^{-1}_1\vecu_{0,1}}{\vecv_{0,1}}_{\Omega_1} -\sum_{i=1}^2 \inner{p_i}{\nabla_i\cdot\vecv_{0,i}}_{\Omega_i} \nonumber\\
    &\qquad= -\sum_{i=1}^2\innerangle{g_{D,i}}{\mathrm{tr}\,\vecv_{0,i}\cdot\vecn_{i}}_{\partial_D\Omega_i}, \qquad\forall\,\vecv_0=[\vecv_{0,1},\vecv_{0,2}]\in {V}_{0}, \label{eq:dualFormVelocity}\\
    &\inner{\mathcal{K}_2^{-1}\left(\vecu_{0,2} + \mathcal{R}_1\lambda_1 + \mathcal{R}_2\lambda_2\right)}{\mathcal{R}_1\nu_1 +\mathcal{R}_2\nu_2}_{\Omega_2} - \inner{p_2}{\nabla_2\cdot\left(\mathcal{R}_1\nu_1 + \mathcal{R}_2\nu_2\right)}_{\Omega_2} \nonumber\\
    &\qquad+\sum_{j=1}^2 \innerangle{\kappa_j^{-1}\lambda_j}{\nu_j}_{\Gamma_j} + \inner{p_1}{\nu_1 + \nu_2 }_{\Omega_1} = 0, \qquad\forall\,\nu=[\nu_1,\nu_2]\in L^2(\Gamma), \label{eq:dualFormMortarFlux}\\
    & \inner{\nabla_2\cdot\left(\vecu_{0,2} + \mathcal{R}_1\lambda_1 + \mathcal{R}_2\lambda_2 \right)}{q_2}_{\Omega_2} + \inner{\nabla_1\cdot\vecu_{0,1}}{q_1}_{\Omega_1} - \inner{\lambda_1 + \lambda_2}{q_1}_{\Omega_1} \nonumber\\
    &\qquad= \sum_{i=1}^2 \inner{f_i}{q_i}_{\Omega_i}, \qquad\forall\,q=[q_1,q_2]\in L^2(\Omega).\label{eq:dualFormPressure}
\end{alignat}
\end{subequations}
\end{definition}

Refer to Appendix~\ref{sec:dualWeakFormProof} for the derivation.

\begin{remark}[Well-posedness] The variational formulation from Definition~\ref{def:dualSingleFrac} can be classified as a saddle point structure, for which well-posedness results have been established for fracture networks, see e.g. Theorem~2.5 from~\cite{boon2018robust}.
\end{remark}

\subsection{A first \textit{a posteriori} error estimate for the primal variable\label{sec:aFirstAposterioriPrimal}}

Having the functional spaces and weak formulations formally introduced, in this section, we provide a first upper bound for an approximation to the primal variable $q = [q_1,q_2]\in H_{0}^1(\Omega) + g$ for the case of a single fracture in the energy norm
\begin{equation}
    \tnorm{q}^2 := \sum_{i=1}^2\, \norm{\mathcal{K}_{i}^\phalf\nabla_i\,q_i}^2_{\Omega_i} + \sum_{j=1}^2\, \norm{\kappa_j^\phalf \left(q_1 - \mathrm{tr}_{\partial_j\Omega_2}\,q_2\right)}_{\Gamma_j}^2. \label{eq:firstEnergyNorm}
\end{equation}

\begin{theorem}[A first upper bound for the primal variable\label{thm:firstUpperBound}] Let $p\in \mixHIzero + g$ be the solution to the primal weak form~\eqref{eq:primalSingleFracture} with $\partial_D\Omega_1$ non-empty. Then for any $q\in \mixHIzero + g$, it holds that
\begin{equation}
    \tnorm{p-q} \leq \sum_{i=1}^2 \eta_{\mathrm{DF},\Omega_i} + \sum_{j=1}^2\eta_{\mathrm{DF},\Gamma_j} + \sum_{i=1}^2 \eta_{\mathrm{R},\Omega_i}, \quad\forall~[\vecv_0,\nu]\in \mixHdiv,\label{eq:firstUpperBound}
\end{equation}
with
\begin{subequations}
\begin{align}
    \eta_{\mathrm{DF},\Omega_1} &= \norm{\mathcal{K}_{1}^{-\frac{1}{2}}\left(\vecv_{0,1}+\mathcal{K}_1\nabla_1\,q_1\right)}_{\Omega_1}, \label{eq:boundDiffusiveFracture} \\
    \eta_{\mathrm{DF},\Omega_2} &= \norm{\mathcal{K}_{2}^{\mhalf}\left(\vecv_{0,2} + \mathcal{R}_1 \nu_1 + \mathcal{R}_2 \nu_2 +\mathcal{K}_2\nabla_2\,q_2\right)}_{\Omega_2}, \label{eq:boundDiffusiveMatrix}\\
    \eta_{\mathrm{DF},\Gamma_1} &= \norm{\kappa_1^{\mhalf}\left(\nu_1 + \kappa_1\left(q_1-\mathrm{tr}_{\partial_1\Omega_2}\,q_2\right)\right)}_{\Gamma_1}, \label{eq:boundDiffusiveInterface1}\\
    \eta_{\mathrm{DF},\Gamma_2} &= \norm{\kappa_2^{\mhalf}\left(\nu_2 + \kappa_2\left(q_1-\mathrm{tr}_{\partial_2\Omega_2}\,q_2\right)\right)}_{\Gamma_2}, \label{eq:boundDiffusiveInterface2}\\
    \eta_{\mathrm{R},\Omega_1} &= 
    C_{\Omega_1}\norm{f_1 - \nabla_1\cdot\vecv_{0,1} + \nu_1 + \nu_2}_{\Omega_1}, \label{eq:boundResidualFracture}\\
    \eta_{\mathrm{R},\Omega_2} &= C_{\Omega_2}\norm{f_2-\nabla_2\cdot\left(\vecv_{0,2}+\mathcal{R}_1\nu_1 + \mathcal{R}_2\nu_2 \right)}_{\Omega_2}, \label{eq:boundResidualMatrix}
\end{align}
\end{subequations}
where
$C_{\Omega_1}$ and $C_{\Omega_2}$ are the permeability-weighted Poincar{\'e}-Friedrichs constants for $\Omega_1$ and $\Omega_2$:
\begin{equation}
C_{\Omega_i} := \sup_{q\in H^1_{0,D}(\Omega_i)}
    \frac{\norm{q}_{\Omega_i}}{\norm{\mathcal{K}_i^{\frac{1}{2}}\nabla_i q}_{\Omega_i}}.
\end{equation}
\end{theorem}

\begin{proof}
Refer to Appendix~\ref{sec:firstThmProof} for the proof.
\end{proof}

\begin{remark}[Nature of the estimators]
The upper bound~\eqref{eq:firstUpperBound} is a guaranteed upper bound for the deviation between the primal solution $p\in \mixHIzero+g$ and an arbitrary approximation $q\in \mixHIzero+g$ in the energy space. There are three types of contributions to the upper bound: (1) diffusive flux estimators \eqref{eq:boundDiffusiveFracture} and \eqref{eq:boundDiffusiveMatrix} measuring the difference between the approximate fluxes $\vecv_{0} + \mathcal{R}\nu \in {V}$ and fluxes obtained from $H^1_0(\Omega)$-potentials $q$, (2) domain coupling estimators \eqref{eq:boundDiffusiveInterface1} and \eqref{eq:boundDiffusiveInterface2} measuring how close the approximate normal fluxes $\nu\in L^2(\Gamma)$ are to the jump in $H^1_0(\Omega)$-potentials $q$, and (3) residual estimators \eqref{eq:boundResidualFracture} and \eqref{eq:boundResidualMatrix} measuring the difference between the exact source term and the divergence of the approximate flux plus the jump in adjacent approximate normal fluxes. An important detail is that the approximate cross-domain fluxes $\nu_1$ and $\nu_2$ enter into the residual estimators of both the higher- and lower-dimensional subdomain. 
\end{remark}

\begin{remark}[Sharpness of the estimates] The estimates above are in principle sharp, as can be shown by standard arguments~\cite{repin2008posteriori}. However, in practice, we will often have access to additional information about the approximate solution (most commonly if it is derived with a local conservation property). This allows for improvements in the residual estimators \eqref{eq:boundResidualMatrix} and \eqref{eq:boundResidualFracture}, as we will show in Section~\ref{sec:evalMajorant}.
\end{remark}

It is clear that even for this fairly simple configuration, the variational formulations (and the analysis in general) can be quite cumbersome. The situation escalates in complexity when intersecting fractures (see Figure~\ref{fig:mixedDecomposition}) are part of the geometric configuration, in particular as the proof of Theorem~\ref{thm:firstUpperBound} relies on all subdomains having some non-vanishing Dirichlet boundary. Indeed, when floating subdomains (e.g., fully embedded fractures or isolated rock domains) are present in the fracture network, the standard procedure used in Theorem~\ref{thm:firstUpperBound} can no longer be applied directly. Thus, in the remainder of the paper, we deal with these challenges in a more general framework. 

\section{Extension to fracture networks\label{sec:extensionToFracNet}}

In this section, we extend the single fracture model to account for several subdomains as part of a general fracture network. Our vocabulary is motivated by the physical case of $n=3$, where the surrounding rock is composed of simply connected 3d subdomains, fractures are simply connected planar 2d subdomains, the intersection between such fractures are 1d lines, and the intersection between fracture intersections are 0d points (see Figure~\ref{fig:mixedDecomposition} for an example with $n=2$).

We start with the classical description and then introduce the mD~notation. The rest of the section is devoted to introducing key tools that are necessary to perform the analysis in an mD~setting.

\subsection{Mixed-dimensional geometric representation}

\begin{figure}[tbp]
    \includegraphics[width=\textwidth]{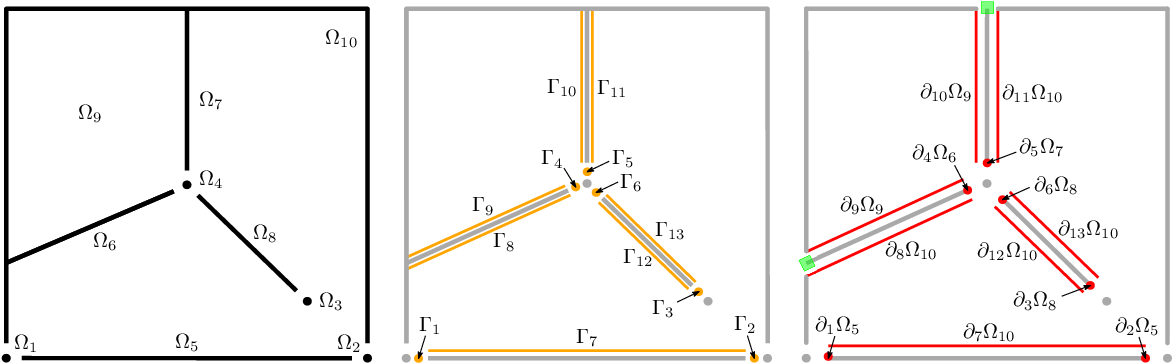}
    \caption{Mixed-dimensional geometric decomposition of a fracture network. Left: The domain $Y$ is decomposed into two 2d matrices ($\Omega_9$ and $\Omega_{10}$), four 1d fractures ($\Omega_5$, $\Omega_6$, $\Omega_7$, and $\Omega_8$), one 0d fracture intersection point ($\Omega_4$), and three 0d fracture end-points ($\Omega_1$, $\Omega_2$, $\Omega_3$). Note that we allow fractures and other lower-dimensional subdomains to form parts of the boundary of the domain (e.g., $\Omega_5$ with its endpoints $\Omega_1$ and $\Omega_2$). Center: Interfaces between subdomains. Right: Subdomain boundaries. Internal boundaries are depicted in red, whereas fracture's boundaries touching the ambient boundary are depicted in green. \label{fig:mixedDecomposition}}
\end{figure}

The derivation of \textit{a posteriori} estimates for generic fracture networks greatly benefits from an mD~decomposition of the domain of interest, and we therefore follow the approach of~\cite{boon2018robust}. We start by considering an $n$--dimensional contractible domain $Y\subset \mathbb{R}^n $, $n\in\{2,3\}$, decomposed into $m$ planar, open and non-intersecting subdomains $\Omega_i$ of different dimensionality $d_i = d(i)$, such that $Y = \cup_{i=1}^m \Omega_i$ (see left Figure~\ref{fig:mixedDecomposition}). The partitioning is constrained such that any $d$-dimensional subdomain (for $d<n$) is always either the intersection of the closure of two or more subdomains of dimension $d+1$, or a cut in a domain of dimension $d+1$. This hierarchical structure excludes e.g., a 1d line or a 0d point embedded directly in a 3d domain. 

We adopt a structure where neighboring subdomains one dimension apart are connected via interfaces, denoted by $\Gamma_j$ for $j\in\{1,\ldots,M\}$. To be precise, let $\Gamma_j$ be the interface between subdomains indexed by $\checkj$ and $\hatj$ of dimension $d$ and $d+1$, respectively. Then $\Gamma_j=\Omega_{\checkj}$ (see center Figure~\ref{fig:mixedDecomposition}), and furthermore, we denote the adjacent boundary of $\Omega_{\hatj}$ by $\Gamma_j = \partial_j \Omega_{\hatj}$. We emphasize that while the internal boundary $\partial_j\Omega_{\hatj}$ is defined to spatially coincide with the interface $\Gamma_j$, which in turn coincides with the lower-dimensional subdomain $\Omega_\checkj$, their distinction is crucial to define variables properly.

To keep track of the connections from subdomains to interfaces, we introduce the sets $\shigh$ and $\slow$, containing the indices of the higher-dimensional (respectively lower-dimensional) neighboring interfaces of $\Omega_i$, as illustrated in the right panel of Figure~\ref{fig:mixedDecomposition}. These sets are dual to $\checkj$ and $\hatj$ defined in the previous paragraph, thus for all $j\in\shigh$, it holds that $\checkj=i$, while for all $j\in\slow$, it holds that $\hatj=i$.

We will be interested in defining functions on the above stated partition of the domain and the interfaces. This motivates us to define the disjoint unions 
\begin{equation}
    \Omega = \bigsqcup_{i=1}^m \Omega_i\qquad \mathrm{and}\qquad \Gamma = \bigsqcup_{j=1}^M \Gamma_j.
\end{equation} 
A complete mixed-dimensional partitioning, including both subdomain and interfaces, is given by $ \Omega \sqcup \Gamma$.

In order to speak of boundary conditions, we introduce the decomposition of the boundary of $\Omega$. Let $\partial\Omega$ be partitioned into its Neumann, Dirichlet, and internal parts. That is, we define $\partial\Omega = \partial_N\Omega \cup \partial_D\Omega \cup \partial_I\Omega$, where $\partial_N\Omega = \cup_{i=1}^m \partial_N\Omega_i$, $\partial_D\Omega = \cup_{i=1}^m \partial_D\Omega_i$, and $\partial_I\Omega = \cup_{i=1}^m \cup_{j\in\slow} \partial_j\Omega_i$.  Finally, to ensure the existence of a unique solution, we require $\partial_D\Omega \neq \emptyset$. 

\subsection{The model problem for a fracture network\label{sec:fractureNetwork}}

Let us now present the model problem valid for $m$ subdomains of dimensionality $0$ to $n$, and $M$ interfaces of dimensionality $0$ to $n-1$. Our model summarizes the derivations given in recent literature \cite{boon2018robust, nordbotten2019unified, keilegavlen2020porepy}. For all domains $\Omega_i$, we consider a scalar pressure $p_i$ together with a flux $\vecu_i$ in the tangent space of the domain. On all interfaces $\Gamma_j$, we consider a scalar coupling flux $\lambda_j$, oriented as positive for flow from the higher dimensional domain $\Omega_{\hatj}$. We will, in this section, assume sufficient regularity that the strong form makes sense, and return to the weak formulation in later sections. The governing equations from the previous section then generalize as
\begin{subequations}
\begin{alignat}{3}
    \hspace{-3mm}\nabla_{i}\cdot\vecu_i - \ssum_{j\in\hat{S}_i}\lambda_{j} &= f_i, \quad&&\mathrm{in}\,\Omega_i, \quad&&i\in \{1,\dots,m\}, \label{eq:conservationGeneral}\\
    \vecu_i &= -\mathcal{K}_i \nabla_{i}\, p_i, \quad&&\mathrm{in}\,\Omega_i, \quad&&i\in \{1,\dots,m\},\,~ d_i\neq 0, \label{eq:darcyGeneral}\\
    \lambda_{j} &= -\kappa_{j} \left(p_\checkj - p_\hatj \right), \quad&&\mathrm{on}\,\Gamma_{j}, \quad&&j\in \{1,\dots,M\}, \label{eq:mortarGeneral} \\
    \vecu_{\hatj}\cdot\vecn_{\hatj} &= \lambda_{j}, \quad&&\mathrm{on}\,\partial_j\Omega_{\hatj}, \quad&&j\in \{1,\dots,M\}, \label{eq:internalBoundaryGeneral}\\
    \vecu_i\cdot\vecn_i &= 0, \quad&&\mathrm{on}\,\partial_N\Omega_i, \quad&&i\in \{1,\dots,m\}, \label{eq:NeumannBoundaryGeneral}\\
     p_i &= g_{D,i}, \quad&&\mathrm{on}\,\partial_D\Omega_i, \quad&&i\in \{1,\dots,m\}. \label{eq:DirichletBoundaryGeneral}
\end{alignat}\label{model:generalClassicalNotation}
\end{subequations}

In~\eqref{eq:conservationGeneral}, the summation captures the contribution of fluxes from the adjacent interfaces to $\Omega_i$, and can be seen as a generalization of the second term in~\eqref{eq:massConservationFracture}. 
Note that for $d_i=n$, the set $\hat{S_i}=\emptyset$, and thus the jump operator, evaluates to zero in the highest-dimensional domains.  Conversely, in \eqref{eq:conservationGeneral}, the differential term $\nabla_i\cdot \vecu_i$ is void whenever $d_i=0$, as there is no tangent space to a point in all subdomains, and indeed, we will not consider the $\vecu_i$ defined on these domains, which justifies why equation \eqref{eq:darcyGeneral} are not applied to 0d domains.

We are now ready to recast the model problem in mD~notation, building on the product space structures introduced in Section~\ref{sec:func_simple}. Let us start by defining the mD~pressure as the ordered collection of subdomain pressures $\gothp := \left[p_i\right]\in C\Omega$, i.e., scalar functions on $\Omega$. 
We now decompose the fluxes as in \eqref{eq:decompose_flux}, so that 
\begin{equation}
\label{eq:decompose}
    \vecu_i = \vecu_{0,i} + \ssum_{j\in\slow}\mathcal{R}_j\lambda_{j}
\end{equation}
such that $\vecu_{0,i}$ satisfies $\vecu_{0,i}\cdot\vecn_i = 0$ for all $j\in \slow$, and where the reconstruction operator is generalized as $\mathcal{R}_j : C\Gamma_j\rightarrow C\Omega_\hatj$ satisfying: 
\begin{equation}
\label{eq:reconstruction2}
    \mathrm{tr}_{\partial_j\Omega_\hatj}\,(\mathcal{R}_j\,\lambda_j) \cdot \vecn_\hatj = 
    \begin{cases}
         \lambda_j & \mathrm{on}\,\partial_j\Omega_\hatj \\
        0 & \mathrm{on}\,\partial\Omega_\hatj\setminus \partial_j\Omega_\hatj
    \end{cases}.
\end{equation}

This allows us to define the mD~flux as the internal (tangential) domain fluxes and (normal) interface fluxes $\gothu := \left[\vecu_{0,i},\lambda_j\right]\in C_0 T\Omega\times C\Gamma$, i.e., the pairing of sections of the tangent bundle $T\Omega$ together with scalar functions on $\Gamma$. By the subscript $C_0 T\Omega$, we indicate that both $\vecu_i\cdot\vecn_i = 0$ on all $\partial_j\Omega_i$, where $j\in \slow$, and also $\vecu_i\cdot\vecn_i = 0$ on $\partial_N\Omega_i$.

We now define a generalized divergence operator $\mixdiv\left(\cdot\right) : C_0 T\Omega\times C\Gamma \rightarrow C\Omega$ which acts on the mD~flux in accordance with the left-hand side of \eqref{eq:conservationGeneral}:
\begin{equation}
    \mixdiv \gothu = \mixdiv \left[\vecu_{0,i}, \lambda_j \right] = \gothq,\label{eq:mD-Divergence}
\end{equation}
where $\gothq = [q_i]\in C\Omega$ is a scalar function for each domain $\Omega_i$, defined by: 
\begin{equation}
q_i := \nabla_i\cdot \left( \vecu_{0,i} + \ssum_{j\in\slow}\mathcal{R}_j\lambda_{j}\right) - \ssum_{j\in\shigh}\lambda_{j}
\end{equation}

Similarly, we define an mD~gradient operator  $\mixgrad\left(\cdot\right) : C\Omega \rightarrow C T\Omega\times C\Gamma$ acting on the mD~pressure in accordance with the right-hand sides of equations \eqref{eq:darcyGeneral} and \eqref{eq:mortarGeneral}: 
\begin{equation}
    \mixgrad \gothp = \mixgrad \left[p_i\right] =\gothv,\label{eq:mD-Gradient}
\end{equation}
where $\gothv = [\vecv_{0,i},\nu_j]\in C T\Omega\times C\Gamma$ has the same structure as the mD~flux (but without the boundary conditions), such that for all $i\in\{1,\ldots,m\}$ and $j\in\{1,\ldots,M\}$, it holds that
\begin{equation}
\label{eq:grad_explicit}
    \nu_j := p_\checkj - p_\hatj, \qquad
    \vecv_{0,i} := \nabla_i\,p_i - \ssum_{j\in\slow}\mathcal{R}_j\nu_{j} .
\end{equation}
Recalling that the full flux $\vecv_i$ is recovered from equation \eqref{eq:decompose}, we note that the second term above is simply the gradient on each subdomain. We will, in Section~\ref{sec:mD_variationalFormulations}, further justify the terminology ``divergence'' and ``gradient'' due to the fact that these operators satisfy an integration-by-parts property with respect to the suitable inner products, and are thus adjoints (subject to appropriate boundary conditions).

Material parameters are collected into the mD~permeability $\mixperm : C T\Omega\times C\Gamma \rightarrow C T\Omega\times C\Gamma$, defined such that for
\begin{equation}
-\mixperm\gothv =  -\mixperm\left[\vecv_{0,i}, \nu_{j}\right] =\gothu, 
\end{equation}
then from the model given in equation \eqref{model:generalClassicalNotation}, we recognize the desired relationships
\begin{equation}
\lambda_j= - \kappa_{j} \nu_{j},
\qquad
\vecu_i = - \mathcal{K}_i 
\vecv_i.
\end{equation}
The second term, corresponding to Darcy's law, can be rewritten in terms of the decomposition $\gothu=[\vecu_{0,i}, \lambda_j]$ from equation \eqref{eq:decompose} as: 
\begin{equation}
\vecu_{0,i} = - \mathcal{K}_i \left(
\vecv_{0,i} + \ssum_{j\in\slow}\mathcal{R}_j\nu_{j}
\right)
- \ssum_{j\in\slow}\mathcal{R}_j\lambda_{j}.
\end{equation}
The presence of the extra terms arising from the decomposition is analogous to that in \eqref{model:generalClassicalNotation}.

We note that the restriction $\gothu \in C_0 T\Omega\times C\Gamma$, implicitly places constraints (depending on the material constants $\mixperm$ and via the definition of $\mixgrad$) on the admissible pressures $\gothp$. This space of admissible pressures can be understood as the domain of the restricted operator $\mixperm\mixgrad : C\Omega \rightarrow C_0 T\Omega\times C\Gamma$. 

In view of the mD~variables and operators defined above, and subject to the right-hand side data $\gothf = [f_i]\in C\Omega$ and the boundary data $\gothg_D = [g_{D,i}]\in C\partial_D\Omega$, a straightforward substitution of definitions shows that problem~\eqref{model:generalClassicalNotation} is equivalent to the concisely stated mD~elliptic problem
\begin{subequations}
    \begin{alignat}{2}
    \gothu &= - \mixperm \mixgrad \gothp, &&\qquad \mathrm{in} \, \Omega\times\Gamma, \label{eq:mD-Darcylaw} \\
    \mixdiv \gothu &= \gothf, &&\qquad \mathrm{in} \, \Omega, \label{eq:mD-massConservation} \\
    \gothp &= \gothg_D, &&\qquad \mathrm{on} \, \partial_D\Omega, \label{eq:mD-DirichletCondtion}
    \end{alignat}
    \label{model:mD-StrongForm}
\end{subequations}
defined for $\gothu\in C_0 T\Omega\times C\Gamma$ and $\gothp\in C\Omega$.

\begin{remark}[Internal Neumann boundaries] For simplicity of exposition, the domain $Y$ is taken as contractible, and $\Omega_i$ is considered a partitioning of $Y$. However, the reader will appreciate that these assumptions can be relaxed. Most importantly, from the perspective of applications (as discussed in Section~\ref{ref:mod_single}), some internal interfaces may be modeled as impermeable, i.e. $\lambda_j = 0$. We refer to the remaining (permeable) interfaces as $\Xi\subset \{0,\ldots,M\}$. The impermeable interfaces are then excluded from the problem, and considered as internal Neumann interfaces. To be precise, we define a reduced disjoint union of interface domains %
\begin{equation*}
    \Gamma=\bigsqcup_{j\in \Xi} \Gamma_j.
\end{equation*}
The internal Neumann boundaries may partition the domain into disconnected parts. We refer to a subdomain as ``Dirichlet-connected'', denoted $i\in \xi$ if either (1) $\partial_D\Omega_i\neq\emptyset$, or (2) there exists some $j\in\shigh$ such that $\hatj\in\xi$, or (3) there exists some $j\in\slow$ such that $\checkj\in\xi$. This allows us to construct a reduced disjoint union of subdomains
\begin{equation*}
\Omega=\bigsqcup_{i\in \xi}\Omega_i.
\end{equation*}
All the derivations in the continuation are equally valid for these reduced product domains.
\end{remark}

\begin{remark}[Extensions to the model equations] The results of this paper can with minor modifications be extended to non-zero Neumann boundary conditions, and with some additional effort to the class of non-planar geometries considered in \cite{boon2020functional}. However, as this generality is typically not needed for applications, we restrict the presentation as indicated above. 
\end{remark}

\subsection{Variational formulations in mixed-dimensional notation\label{sec:mD_variationalFormulations}}

Before writing the variational formulations in mD~notation, let us first define the relevant mD~inner products and norms. Consider the following inner-products
\begin{alignat}{2}
    &\inner{\gothq}{\gothr}_\Omega = \sum_{i=1}^m \inner{q_i}{r_i}_{\Omega_i} \qquad\forall~\gothq=[q_i], \gothr=[r_i]\in L^2\Omega,  \label{eq:innerSubdomain}\\[1.5mm]
    &\inner{\gothv}{\gothw}_{\Omega,\Gamma} = \sum_{i=1}^m \Bigg( \inner{\left(
\vecv_{0,i} + \ssum_{j\in\slow}\mathcal{R}_j\nu_j\right)}{\left(
\vecw_{0,i} + \ssum_{j\in\slow}\mathcal{R}_j\mu_{j}
\right)}_{\Omega_i} \nonumber \\
&~~ + \sum_{j\in\check{S}_i}\innerangle{\nu_j}{\mu_j}_{\Gamma_j}\Bigg) \qquad\forall~\gothv=[\vecv_{0,i},\nu_j],\gothw=[\vecw_{0,i},\mu_j]\in L^2 T\Omega \times L^2\Gamma, \label{eq:innerForest}\\[1.5mm]
    &\innerangle{\gothq}{\gothr}_{\partial_X\Omega} = \sum_{i=1}^m \innerangle{q_i}{r_i}_{\partial_X\Omega_i}  \quad\forall~\gothq=[q_i],\gothr=[r_i]\in L^2 \partial_X\Omega, \label{eq:innerDirichlet}
\end{alignat}
and their respective induced norms~
\begin{equation}
    \norm{\gothq}^2_\Omega = \inner{\gothq}{\gothq}_\Omega, \qquad \norm{\gothv}_{\Omega,\Gamma}^2 = \inner{\gothv}{\gothv}_{\Omega,\Gamma}, \qquad \norm{\gothq}_{\partial_X\Omega}^2 = \innerangle{\gothq}{\gothq}_{\partial_X\Omega}.
\end{equation}

With these inner products, the previously defined mD~divergence satisfy the following integration-by-parts formula \cite{boon2018robust,boon2020functional} whenever $\gothv\in C T\Omega\times C\Gamma$ and $\gothq\in C\Omega$.
\begin{equation}
\label{eq:int-by-parts}
    \inner{\gothq}{\mixdiv\gothv}_\Omega + \inner{\mixgrad\gothq}{\gothv}_{\Omega,\Gamma}
    =
    \innerangle{\mixTraceD\gothq}{\mixTraceD\gothv}_{\partial_D\Omega} + \innerangle{\mixTraceN\gothq}{\mixTraceN\gothv}_{\partial_N\Omega}.
\end{equation}
In the above the restriction to the boundary is denoted $\mixTraceX({\cdot})$ (for $X=D,N$), which depending on context acts as the boundary values of pressure variables, $\mixTraceX({\cdot}) : C\Omega \rightarrow C\partial_X\Omega$, or the normal component of flux variables, $\mixTraceX({\cdot}) : C T\Omega\times C\Gamma \rightarrow C\partial_X\Omega$.

From the product structure in the definition of the $C$ and $L^2$ spaces, the continuous spaces inherit their density from the individual subdomains to the product spaces on $\Omega$ and $\Gamma$. We can thus follow standard procedures to obtain weak extensions of the mD~differential operators, the boundary restriction (trace) operators, and the corresponding function spaces~\cite{adams2003sobolev, pedersen1989analysis, arnold2018FEEC}. We elaborate this below.

Due to the density of $C_0 T\Omega \times C\Gamma$ in $L^2 T\Omega \times L^2\Gamma$, the mD~divergence from Section~\ref{sec:fractureNetwork} is a densely defined unbounded linear operator on the latter space $\mixdiv : L^2\Omega \rightarrow L^2 T\Omega\times L^2\Gamma$. Let us now (temporarily) use the notation $(T,\mathrm{dom}(T))$ to emphasize that an operator $T$ has domain of definition $\mathrm{dom}(T)$, and we denote the adjoint operator with respect to the $L^2$ inner product by an asterisk. 

We recall that the Neumann boundary is incorporated into the definition of the continuous flux spaces $C_0 T\Omega\times C\Gamma$, thus the last term in the integration-by-parts formula \eqref{eq:int-by-parts}, is zero. Hence, we can define a weak mD~gradient and the corresponding space of weakly mD~differentiable functions with zero trace on the Dirichlet boundary $H^1_0$ by considering the adjoint: 
\begin{equation}
    (\mixgrad, H^1_0(\Omega)) := (\mixdiv, C_0 T\Omega\times C\Gamma)^*.
\end{equation}
Clearly, $C_0\Omega \subseteq H^1_0(\Omega)$, and thus it is appropriate to consider $(\mixgrad, H^1_0(\Omega))$ as a weak gradient. Moreover, the domain of definition simply corresponds to the standard $H^1_0(\Omega_i)$
on each domain, where the subscript zero indicates zero trace on all Dirichlet boundaries. Thus $H^1_0(\Omega) = \prod_{i=1}^m H^1_0(\Omega_i)$, which generalizes \eqref{eq:prod2}. 

Considering the integration-by-parts formula again, the weak mD~divergence and the corresponding space of flux functions with divergence in $L^2$ and zero trace on the Neumann boundary $H(\mathrm{div}; \Omega,\Gamma)$ can be defined as 
\begin{equation}
    (\mixdiv, H(\mathrm{div}; \Omega,\Gamma)) := (\mixgrad, H^1_0)^*.
\end{equation}
Again $C_0 T\Omega\times C\Gamma\subseteq H(\mathrm{div}; \Omega,\Gamma)$, and it is appropriate to consider \\$(\mixdiv, H(\mathrm{div}; \Omega,\Gamma))$ as a weak divergence. 
This domain of definition of the weak divergence has the interpretation of $H_0(\mathrm{div}; \Omega_i)$ on all subdomains $\Omega_i$ (where the subscript zero indicates zero trace on all boundaries except for Dirichlet boundaries), and $L^2(\Gamma_j)$ spaces on all interfaces $\Gamma_j$. Thus $H(\mathrm{div}; \Omega,\Gamma) = \prod_{i=1}^m H_0(\mathrm{div};\Omega_i) \times \prod_{i=1}^M L^2(\Gamma_j)$, which generalizes \eqref{eq:firstHdiv}.

Due to the above identification of $H^1(\Omega)$ and $H(\mathrm{div}; \Omega,\Gamma)$ in terms of product spaces of standard function spaces on subdomains, we extend the definition of the boundary restriction operators $\mixTraceX({\cdot})$ to trace operators on the weak spaces by requiring that they coincide with the standard trace operators on subdomains. 

In the continuation, we will always consider the weak mD~gradient and divergence, and denote these simply by $\mixgrad$ and $\mixdiv$, respectively. Similarly, we will always consider the boundary restrictions as trace operators. The above definitions of weak mD~gradient and divergence operators, and their adjoint property on the above weak spaces, has the following statements of the primal and dual weak formulations of equations \eqref{model:mD-StrongForm} as a direct consequence: 

\begin{definition}[Mixed-dimensional primal weak form\label{def:mD_primal}] Let $\gothg\in H^1(\Omega)$. Then find $\gothp \in H_0^1(\Omega) + \gothg$ such that
\begin{equation}
    \inner{\mixperm\mixgrad\gothp}{\mixgrad\gothq}_{\Omega,\Gamma} = \inner{\gothf}{\gothq}_{\Omega} \qquad \forall \, \gothq \in H_0^1(\Omega). \label{eq:mD_weakPrimal}
\end{equation}
\end{definition}

\begin{definition}[Mixed-dimensional dual weak form\label{def:mD_dual}] Find $(\gothu,\gothp) \in \mixHdiv\times L^2(\Omega)$ such that
\begin{subequations}
    \begin{alignat}{2}
        \inner{\mixperm^{-1} \gothu}{\gothv}_{\Omega,\Gamma} - \inner{\gothp}{\mixdiv\gothv}_{\Omega} &= \innerangle{\gothg_D}{\mixTraceD{\gothv}}_{\partial_D\Omega} \qquad&&\forall\,\gothv\in\mixHdiv, \label{eq:mD_weakDual1}\\
        \inner{\mixdiv \gothu}{\gothq}_{\Omega} &= \inner{\gothf}{\gothq}_{\Omega}  \qquad&&\forall\,\gothq\in L^2(\Omega). \label{eq:mD_weakDual2}
    \end{alignat}
\label{eq:mD_weakDual}
\end{subequations}
\end{definition}

The above weak forms of the mixed-dimensional elliptic problem are well-posed for bounded coefficients \cite{boon2020functional}, in the sense that there exist positive constants $\mixperm_0$ and $\mixperm_\infty$ such that:
\begin{equation}
    \sup_{\gothv\in\mixHdiv} \frac
    {\inner{\mixperm \gothv}{\gothv}_{\Omega,\Gamma}}
    {\mixperm_\infty\norm{\gothv}^2_{\Omega,\Gamma}} \leq 1\leq
    \inf_{\gothv\in\mixHdiv} \frac
    {\inner{\mixperm \gothv}{\gothv}_{\Omega,\Gamma}}
    {\mixperm_0\norm{\gothv}^2_{\Omega,\Gamma}}.
\end{equation}
The solutions of the primal and dual weak formulations are equivalent, and define true solutions $\gothp\in\mixHIzero+\gothg$ and $\gothu\in\mixHdiv$ against which the approximate solutions will be measured in later sections.  

\section{Functional analysis tools}
\label{sec:functools}
In this section, we summarize the main functional analysis tools we will exploit for the \textit{a posteriori} analysis. 

\subsection{Poincar{\'e}-type inequalities\label{sec:Poincare&Weighted}}

We recall the following weighted Poincar{\'e} inequalities:

\begin{lemma}[Permeability-weighted Poincar{\'e}-Friedrichs inequalities\label{lemma:Poincare}] There exist constants $C_\Omega \geq C_{\Omega_i} \geq C_K$ such that
\begin{subequations}
\begin{alignat}{3}
    \norm{\gothq}_{\Omega,\Gamma} &\leq C_{\Omega,\Gamma} ~\norm{\mixperm^{\phalf}\mixgrad\gothq}_{\Omega,\Gamma} \qquad&&\forall~\gothq\in\mixHIzero, \label{eq:CP_mixdim}\\
    \norm{q}_{\Omega_i} &\leq C_{\Omega_i}~\norm{\mathcal{K}_i^{\phalf}\nabla_{i} q}_{\Omega_i} \qquad&&\forall~q\in H^1_0(\Omega_i),\quad && \mathrm{if }\, \partial_D\Omega_i \neq \emptyset,
    \label{eq:CP_subdomain}\\
    \norm{q -\tilde{q}_{\Omega_i}}_{\Omega_i} &\leq {C}_{\Omega_i}~\norm{\mathcal{K}_i^{\phalf}\nabla_{i} q}_{\Omega_i} \qquad&&\forall\,q\in H^1(\Omega_i),\quad && \mathrm{if }\, \partial_D\Omega_i = \emptyset, \label{eq:CP_subdomainZeroMean}\\
    \norm{q - \tilde{q}_K}_{K} &\leq {C}_{K}~\norm{\mathcal{K}_i^{\phalf}\nabla_{i} q}_{K} \qquad&&\forall~q\in H^1(K),\quad && \mathrm{where}\, K\subset\Omega_i. \label{eq:CP_local}
\end{alignat}
\end{subequations}
Here, we denote by $\tilde{q}_{\Omega_i}$ and $\tilde{q}_{K}$ the mean value of $q$ over the subdomain $\Omega_i$ and an arbitrary $d_i$-simplex $K\subset\Omega_i$, respectively.
\end{lemma}

We refer to $C_{\Omega,\Gamma}$ as the mixed-dimensional permeability-weighted Poincar{\'e}-Friedrichs constant (whose existence was shown in \cite{boon2020functional}), $C_{\Omega_i}$ is the standard subdomain permeability-weighted Poincar{\'e}-Friedrichs constant, and ${C}_{K}$ is a local permeability-weighted Poincar{\'e}-Friedrichs constant.

It is important to mention that concrete values of $C_{\Omega_i}$ are available only for a limited set of geometries, see e.g.,~\cite{carstensen2000constants, repin2012computable, vohralik2005Poincare}. An upper bound exists for convex domains, and thus for a simplex $K\subset\Omega_i$ we have~\cite{payne1960optimal, bebendorf2003poincare}
\begin{equation}
    C_K \leq \frac{\mathrm{diam}(K)}{\pi c_K}
\end{equation}
where $c_K$ is the lower bound on the permeability within $K$: 
\begin{equation}
    c_K = \inf_{\substack{\vecx\in K\\ \vecv\in TK_\vecx}} \frac{(\mathcal{K}_i(\vecx) \vecv)\cdot \vecv}{\norm{\vecv}}^2
\end{equation}

The importance of this is understood if $K$ is an element of a simplicial partition of $\Omega_i$, in which case $C_{K}$ scales with the mesh size $h_K=\mathrm{diam}(K)$. This allows for super-convergent properties of residual estimators for some locally mass-conservative approximations~\cite{vohralik2010unified, ern2015polynomial, ernGuaranteed2017}. We analyze these cases with further details in Section~\ref{sec:evalMajorant} and Remark~\ref{rem:superconvergence}.

\subsection{Conforming flux spaces}

It is often possible to verify that an approximate solution $\gothv\in\mixHdiv$ satisfies some stronger conservation property, that is to say, that there is some space $U\subseteq L^2$ such that 
\begin{equation}
    \mixdiv\gothv - \gothf\in U
\end{equation}
This allows for the construction of stronger \textit{a posteriori} estimates, and as such, we formalize this concept as a generalization of $\mixHdiv$ to ``$U$-conforming flux spaces'':
\begin{definition}[Conforming mD~flux space] Let $\mixHdivU \subset \mixHdiv$ be a $U$-conforming flux space, in the sense of
\begin{equation}
    \mixHdivU = \left\lbrace \gothv\in\mixHdiv : \gothf - \mixdiv\gothv \in U\right\rbrace.
\end{equation}
\end{definition}

To exploit the conforming flux spaces, we must construct certain projected $\mixHI$ spaces. Consider therefore $U$ as some subspace of $L^2(\Omega)$ and define $U^\perp$ to be its orthogonal complement:
\begin{equation}
    U^\perp := \lbrace \gothq \in L^2(\Omega) : \inner{\gothq}{\gothr}_\Omega  = 0 ~~\forall~ \gothr\in U\rbrace.
\end{equation}
Moreover, let $\pi_{U^\perp}$ be the $L^2$--projection onto $U^\perp$, such that for any $\gothr\in L^2(\Omega)$, $\pi_{U^\perp}\gothr\in U^\perp$ satisfies the orthogonality property:
\begin{equation}
    \inner{\gothr-\pi_{U^\perp}\gothr}{\gothq}_{\Omega} = 0 \quad\forall~\gothq\in U^\perp.\label{eq:ortho}
\end{equation}

Consider now the projected $\mixHIzero$ space denoted $W\subset L^2(\Omega)$, defined as the range of $\pi_{W} := (I - \pi_{U^\perp}) : \mixHIzero \to L^2(\Omega)$, and let the norm of $W$ be defined as a weighted $L^2$-norm with nonnegative weights $\mu \in L^\infty(\Omega)$
\begin{equation}
    \norm{\gothq}_{W,\mu} := \norm{\mu \gothq}_\Omega \qquad\forall~\gothq\in W,\label{eq:weightedNorm}
\end{equation}
which are defined within the class $\calC_W$ with unit Poincar{\'e} constants:
\begin{equation}
    \calC_W = \left\lbrace \mu \in L^{\infty}(\Omega) : \sup_{\gothq\in\mixHIzero} \frac{\norm{\pi_{W} \gothq}_{W,\mu}}{\norm{\mixperm^{\phalf}\mixgrad \gothq}_{\Omega,\Gamma}} \leq 1 \right\rbrace.
\end{equation}
Indeed, such classes exist in the literature of Poincar{\'e} inequalities for weighted norms, see e.g.,~\cite{rathmair2019poincare, pechstein2013weighted}. Note that a trivial member of $\calC_W$ is the inverse of the permeability-weighted mD~Poincar{\'e}-Friedrichs constant $\mu(\vecx)=C_{\Omega,\Gamma}^{-1}$. As we will see in Sections~\ref{sec:GAE} and~\ref{sec:evalMajorant}, the concrete choice of the space $U$ and the corresponding weights $\mu$ will directly impact the strength of the estimates.

\begin{remark}[On the space $\mixHdivU$] The conforming mD~flux spaces allow us to obtain sharper estimates in Section~\ref{sec:posteriori}. However, it is important to note that the standard case $U=L^2(\Omega)$ is included in our definition, for which the orthogonal complement is void, and the projection $\pi_W =I$; thus $W = \mixHIzero$. This and other cases are elaborated in more detail in Sections~\ref{sec:NC} to \ref{sec:EC}.
\end{remark}

\subsection{Bilinear forms and energy norms\label{sec:energyNorms}}

For the \textit{a posterior}i analysis, we will need the next two mD~bilinear forms and their induced energy norms
\begin{alignat}{3}
    \mathfrak{B}(\gothq, \gothr) &= \inner{\mixperm \mixgrad \gothq}{\mixgrad \gothr}_{\Omega, \Gamma}, &&\quad\tnorm{\gothq}^2 = \mathfrak{B}(\gothq,\gothq) = \norm{\mixperm^{\phalf}\mixgrad\gothq}_{\Omega, \Gamma}^2 \quad&&\forall~\gothq,\gothr \in H_0^1(\Omega), \label{eq:bilinearB}\\
    \mathfrak{A}(\gothv, \gothw) &= \inner{\gothv}{\mixperm^{-1}\gothw}_{\Omega, \Gamma}, &&\quad\tnormstar{\gothv}^2 = \mathfrak{A}(\gothv,\gothv) = \norm{\mixperm^{\mhalf}\gothv}_{\Omega, \Gamma}^2 \quad&&\forall~\gothv,\gothw \in L^2T\Omega\times L^2\Gamma,\label{eq:bilinearA}
\end{alignat}
which are related via
\begin{equation}
    \tnorm{\gothq} = \tnormstar{\mixperm\mixgrad\gothq} \qquad \forall\, \gothq\in H_0^1(\Omega).
    \label{eq:normsRelationship}
\end{equation}

We also define the \textit{full} norm for a mixed-dimensional pair of primal and dual variables as
\begin{equation}
\fullnorm{\gothq,\gothv} := \tnorm{\gothq} + \tnormstar{\gothv} + \norm{\mu^{-1} \mixdiv\gothv}_{\Omega} \quad\forall\, (\gothq,\gothv)\in \mixHIzero\times\mixHdivU.
\label{eq:combinedNorm}
\end{equation}
Note that the last norm will depend on the eventual choice of $\mu^{-1}$, which we emphasize must be from the class $\mu \in \calC_W$, as defined in the preceding section. 

\section{A posteriori error estimates\label{sec:posteriori}}

This section is devoted to obtaining the error bounds for our model problem. First, we provide general abstract estimates, and later 
we focus on the evaluation of the different bounds.

\subsection{General abstract estimates\label{sec:GAE}}

Let us now present the general abstract bounds. We formalize the main results presented in Section~\ref{sec:extensionToFracNet} and extend the ones presented in Theorem~\ref{thm:firstUpperBound} in the following theorem.

\begin{theorem}[General abstract a posteriori error bounds\label{thm:Abstract}] Let the error majorant be defined as%
\begin{align}
    \mathcal{M}(\gothq,\gothv,\gothf,\mu) &:= \eta_{\mathrm{DF}}(\gothq,\gothv) + \eta_{\mathrm{R}}(\gothv,\gothf,\mu) ,
\end{align}
where
\begin{equation}
    \eta_{\mathrm{DF}}(\gothq,\gothv) := \tnormstar{\gothv+\mixperm\mixgrad\gothq} \quad \mathrm{and} \quad 
    \eta_{\mathrm{R}}(\gothv,\gothf,\mu) := \norm{\mu^{-1}(\gothf-\mixdiv\gothv)}_{\Omega},
\end{equation}
valid for all $\gothq \in \mixHIzero + \gothg$ and $\gothv\in\mixHdivU$. Then, the following a posteriori error estimates hold.

\noindent (1) Let $\gothp \in \mixHIzero + \gothg$ be the solution to~\eqref{eq:mD_weakPrimal} and $\gothq\in\mixHIzero + \gothg$ be arbitrary. Then
\begin{equation}
    \tnorm{\gothp - \gothq} \leq \mathcal{M}_{\gothp}^{\oplus} =  \mathcal{M}(\gothq,\gothv,\gothf,\mu) \qquad \forall~\gothv\in\mixHdivU, \label{eq:primalAbstract}
\end{equation}
where $\mathcal{M}_\gothp^{\oplus}$ is the upper bound of the error for the primal variable.

\noindent (2) Let $\gothu\in\mixHdiv$ be the solution to~\eqref{eq:mD_weakDual} and $\gothv\in\mixHdivU$ be arbitrary. Then
\begin{equation}
    \tnormstar{\gothu - \gothv} \leq \mathcal{M}_{\gothu}^{\oplus} = \mathcal{M}(\gothq,\gothv,\gothf,\mu) \qquad \forall~\gothq\in \mixHIzero+\gothg, \label{eq:dualAbstract}
\end{equation}
where $\mathcal{M}_{\gothu}^{\oplus}$ is the upper bound of the error for the dual variable.

\noindent (3) Let $\gothp \in \mixHIzero+\gothg$ be the solution to~\eqref{eq:mD_weakPrimal} and $\gothu\in\mixHdiv$ be the solution to~\eqref{eq:mD_weakDual}, and let $(\gothq,\gothv) \in (\mixHIzero+\gothg) \times \mixHdivU$ be arbitrary. Then,
\begin{equation}
    \mathcal{M}(\gothq,\gothv,\gothf,\mu) = \mathcal{M}^\ominus_{\gothp,\gothu} \leq \fullnorm{\gothp-\gothq,\gothu-\gothv} \leq \mathcal{M}_{\gothp,\gothu}^\oplus = 2 \mathcal{M}(\gothq,\gothv,\gothf,\mu) + \eta_{\mathrm{R}}(\gothv,\gothf,\mu), \label{eq:combinedAbstract}
\end{equation}
where $\mathcal{M}^\ominus_{\gothp,\gothu}$ and $\mathcal{M}^\oplus_{\gothp,\gothu}$ are the lower and upper bounds of the error for the primal-dual variable.
\end{theorem}

\begin{proof} Due to the construction of mixed-dimensional product spaces and the adjoint property of the weak differential operators, the proof from the mono-dimensional case can (to a large extent) be applied directly~\cite{repin2007mixed}. A notable deviation from the standard proofs is the use of conforming flux spaces, and the inclusion of the Poincar{\'e}-constants in the weights $\calC_W$. The full proof is included for completeness in Appendix~\ref{sec:thmProof}.
\end{proof}

\begin{remark}[Non-conforming approximations] Referring again to the general setting of mD~calculus, it has been shown that the differential operators form part of a cochain complex, and that an mD~Helmholtz decomposition exists~\cite{boon2020functional}. Thus, by realizing the above constructions as Hilbert complexes, the above error bounds can be extended also to non-conforming approximations following, e.g., Theorem~4.7 of~\cite{pauly2020solution}. However, as a main objective of our work is to obtain bounds based on conforming properties of the approximations, we will not pursue non-conforming approximations in this work.  
\end{remark}

\subsection{Evaluation of the majorant\label{sec:evalMajorant}}

The aim of this section is to provide concrete forms of the majorant $\mathcal{M}(\gothq,\gothv,\gothf,\mu)$ from Theorem~\ref{thm:Abstract} depending upon the choices of the weights $\mu$. For this purpose, consider once again the definition of the majorant
\begin{align}
    &\mathcal{M}(\gothq,\gothv,\gothf,\mu) = \eta_{\mathrm{DF}}(\gothq,\gothv) + \eta_{\mathrm{R}}(\gothv,\gothf,\mu) \nonumber \\
    &\qquad\forall~\gothq=[q_i]\in H_0^1(\Omega)+\gothg,~\gothv=[\vecv_{0,i},\nu_{j}]\in\mixHdivU.
\end{align}

The estimation of the first term $\eta_{\mathrm{DF}}(\gothq,\gothv)$ is independent of the weights $\mu$. Indeed, by applying \eqref{eq:bilinearA}, it is straightforward to see that 
\begin{align}
\label{eq:df_loc}
 \eta_{\mathrm{DF}}^2(\gothq,\gothv) &= \sum_{i=1}^m \left(\sum_{K\in\Tau_{\Omega_i}}\norm{\calK^{\mhalf}_i \Bigg(
\vecv_{0,i} +  \sum_{j\in\slow}\mathcal{R}_j\nu_{j}
\Bigg) + \calK^{\phalf}_i \nabla_{i}q_i}_{K}^2 \right. \nonumber \\
& \left. \hspace{5mm}+ \sum_{j\in\slow} \sum_{K\in\Tau_{\Gamma_{j}}} \norm{\kappa_{j}^{\mhalf} \nu_{j} + \kappa_{j}^{\phalf}\left(q_{\checkj} - \trace q_{\hatj}\right)}^2_{K}\right) \nonumber \\
 &= \sum_{i=1}^m \left(\sum_{K\in\Tau_{\Omega_i}} \eta^2_{\mathrm{DF}_\parallel,K} + \sum_{j\in\slow} \sum_{K\in\Tau_{\Gamma_{j}}} \eta^2_{\mathrm{DF}_\perp,K} \right).
\end{align}
The terms $\eta_{\mathrm{DF}_\parallel,K}$ and $\eta_{\mathrm{DF}_\perp,K}$ measure the diffusive flux error in the tangential and normal directions associated with the subdomain element $K\in\Tau_{\Omega_i}$ and the mortar element $K\in\Tau_{\Gamma_{j}}$, respectively.

To complete the evaluation of the majorant, we are left with the estimation of $\eta_{\mathrm{R}}(\gothv,\gothf,\mu)$, which depends on the choices of $\mu$. Recall that this term measures the mismatch in satisfying the conservation equation in each subdomain. To be precise, there are four main types of conforming fluxes; Standard $L^2$-conforming,  subdomain conservation, grid level (local) conservation, and point-wise. The quality of the residual balance can be verified explicitly \textit{before} applying the \textit{a posteriori} estimates, and thus is not considered an assumption in the theory. Below, we make precise the aforementioned cases.

\subsubsection{No mass-conservation\label{sec:NC}}
\label{sec:nocon}
Assume nothing is known about the approximation of the residual terms beyond the $L^2$ structure. We indicate this case by the abbreviation ``NC'', and set $U_{\mathrm{NC}}=L^2$, and $\gothv \in \mixHdivUa{\mathrm{NC}}=\mixHdiv$. Then $U_{\mathrm{NC}}^\perp = 0$, which implies that $\pi_{W} = I$, and $W = \mixHIzero$. Then, \textit{a priori}, we only know the global (mixed-dimensional) Poincar{\'e}~\eqref{eq:CP_mixdim}, i.e., we have no better weight than setting ${\mu}(\vecx) = C_{\Omega, \Gamma}^{-1}$ for $\vecx\in\Omega$.

Using~\eqref{eq:bilinearB} and the mD~Poincar{\'e} inequality~\eqref{eq:CP_mixdim}, one obtains the following bound, which is the weakest bound available within the class of bounds considered in this paper:
\begin{align}
    \eta_{\mathrm{R}}^2 &\leq C^2_{\Omega, \Gamma} \sum_{i=1}^m \sum_{K\in\Tau_{\Omega_i}} \norm{f_i - \nabla_i\cdot\Bigg(
\vecv_{0,i} + \sum_{j\in\slow}\mathcal{R}_j\nu_{j}
\Bigg)
+ \sum_{j\in\shigh}\nu_{j}}_K^2 \nonumber 
\\
&=  \sum_{i=1}^m \sum_{K\in\Tau_{\Omega_i}} \eta_{\mathrm{R},K;\mathrm{NC}}^{2} = \eta_{\mathrm{R};\mathrm{NC}}^2,
\end{align}
Here, $\eta_{\mathrm{R},K;\mathrm{NC}}$ denotes the local residual error for non-conservative approximations. The majorant when mass conservation cannot be guaranteed at any level is then given by,
\begin{equation}
    \mathcal{M}_\mathrm{NC}(\gothq,\gothv,\gothf) = \eta_{\mathrm{DF}}(\gothq,\gothv) + \eta_{\mathrm{R};\mathrm{NC}}(\gothv,\gothf),\label{eq:M_NC}
\end{equation}
and it follows from the above that this is an upper bound, $\mathcal{M}\leq \mathcal{M}_\mathrm{NC}$.

\subsubsection{Subdomain mass-conservation}
\label{sec:submasscon}

Due to the structure of the equations, where interface fluxes are stated explicitly, many approximations will have mass conservation satisfied in a subdomain level, which is in a sense a compatibility condition on the floating domains $\Omega_i$. We indicate this case by the abbreviation ``SC''. In particular, the divergence $\gothr=[r_i]=\mixdiv \gothv \in U_{\mathrm{SC}}$ satisfies for all $i\in\{1,\ldots,m\}$ where $\partial_D\Omega_i =\emptyset$, 
\begin{equation}
    \inner{r_i}{1}_{\Omega_i} = \inner{f_i}{1}_{\Omega_i}. \label{eq:subdomainMassConservation}
\end{equation}
Thus, by definition $U_{\mathrm{SC}}^\perp$ is the space of constants over the floating subdomains $\Omega_i$, and the space $W$ is the space of $\mathring{H}^1(\Omega_i)$ functions, with zero mean if $\partial_D\Omega_i = \emptyset$. 

This case represents an improvement relative to the previous one, in the sense that we can now employ the subdomain Poincar{\'e} constants instead of the mD~constant. Let us make this point precise in the following lemma.

\begin{lemma}\label{lemma:localPoincare} Let $W = \prod_{i=1}^m \mathring{H}^1(\Omega_i)$, where
\begin{equation}
    \mathring{H}^1(\Omega_i) = \left\lbrace q_i\in H^1_0(\Omega_i)\, |\, \inner{q_i}{1}_{\Omega_i} = 0\, \mathrm{if}\, \partial_D\Omega_i = \emptyset \right\rbrace.\label{eq:zeroMean}
\end{equation}
Then, $\mu(\vecx) = C_{\Omega_i}^{-1}$ for $\vecx \in \Omega_i$ belongs to the class $\mathcal{C}_W$, where $C_{\Omega_i}$ is the permeability-weighted Poincar{\'e}-Friedrichs constants defined in Lemma \ref{lemma:Poincare}.
\end{lemma}
\begin{proof} Using the Poincar{\'e} inequality \eqref{eq:CP_subdomainZeroMean} and the fact that the sum of broken norms is weaker than the full norm, the following result holds  
\begin{align*}
    &\sup_{\gothq\in\mixHIzero} \frac{\norm{\pi_{W}\gothq}_{W,\mu}}{\norm{\mixperm^{\phalf}\mixgrad\gothq}_{\Omega, \Gamma}} = \sup_{\substack{\gothq\in \mixHIzero \\ \norm{\mixperm^{\phalf}\mixgrad\gothq}_{\Omega, \Gamma}=1}} \norm{\pi_{W_\Omega}\gothq}_{W_\Omega,\mu}\\
    &= \sup_{\substack{\gothq\in\mixHIzero \\ \norm{\mixperm^{\phalf}\mixgrad\gothq}_{\Omega, \Gamma}=1}}~ \left(
    \sum_{\substack{i=1 \\ \partial_D\Omega_i \neq \emptyset}}^m \norm{C^{-1}_{{\Omega_i}} q_i}_{\Omega_i}
    +
    \sum_{\substack{i=1 \\ \partial_D\Omega_i = \emptyset}}^m \norm{C^{-1}_{{\Omega_i}} \left(q_i- \frac{1}{\abs{\Omega_i}} \inner{q_i}{1}_{\Omega_i}\right)}_{\Omega_i}
    \right) \\ 
    &
    \leq \sup_{\substack{\gothq\in \mixHIzero \\ \norm{\mixperm^{\phalf}\mixgrad\gothq}_{\Omega, \Gamma}=1}} \sum_{i=1}^m \norm{\mathcal{K}^{\phalf}\nabla_{i} q_i}_{\Omega_i} \leq 1.
\end{align*}
\end{proof}

In view of Lemma~\ref{lemma:localPoincare}, $\eta_{\mathrm{R}}$ can be bounded as
\begin{align}
    \eta_{\mathrm{R}}^2 &\leq  \sum_{i=1}^m C^2_{\Omega_i} \sum_{K\in\Tau_{\Omega_i}} \norm{f_i - \nabla_i\cdot\Bigg(
\vecv_{0,i} + \sum_{j\in\slow}\mathcal{R}_j\nu_{j}
\Bigg)
+ \sum_{j\in\shigh}\nu_{j}}_K^2 \nonumber \\
&= \sum_{i=1}^m \sum_{K\in\Tau_{\Omega_i}} \eta_{\mathrm{R},K;\mathrm{SC}}^2 = \eta_{\mathrm{R};\mathrm{SC}}^2,
\end{align}
where $\eta_{\mathrm{R},K;\mathrm{SC}}$ are the local residual estimators for subdomain mass-conservative approximations. 
The majorant for this case is given by
\begin{equation}
    \mathcal{M}_\mathrm{SC}(\gothq,\gothv,\gothf) = \eta_{\mathrm{DF}}(\gothq,\gothv) + \eta_{\mathrm{R};\mathrm{SC}}(\gothv,\gothf).\label{eq:M_SC}
\end{equation}
This estimate is sharper than that in the preceding section, since $C_{\Omega_i}\leq C_{\Omega,\Gamma}$, thus whenever the assumptions of this section are satisfied, it holds that $\mathcal{M}\leq\mathcal{M}_\mathrm{SC}\leq\mathcal{M}_\mathrm{NC}$.

Note that~\eqref{eq:M_SC} is identical in structure to the residual estimators~\eqref{eq:boundResidualFracture} and \eqref{eq:boundResidualMatrix} obtained in Theorem~\ref{thm:firstUpperBound}. However, they are fundamentally different in the sense that \eqref{eq:M_SC} do not require all subdomains to have a non-empty Dirichlet part but rather mass to be conserved in each subdomain $\Omega_i$.

\subsubsection{Local mass-conservation}
\label{sec:LMC}

By choice of numerical method, it is often easy to verify that mass is conserved on an element basis in a subdomain partition. We indicate this case by the abbreviation ``LC''. As in the preceding section, this implies that the divergence $\gothr=[r_i]=\mixdiv \gothv \in U_{\mathrm{LC}}$ then satisfies for all $K\subset \Tau_{\Omega_i}$ that
\begin{equation}
    \inner{r_i}{1}_{K} = \inner{f_i}{1}_{K}, \label{eq:localMassConservation}
\end{equation}
where $\Tau_{\Omega_i}$ denotes a finite partition of $\Omega_i$ (typically a simplicial grid).
In this case, $U_{\mathrm{LC}}$ contain functions having zero mean on each element $K\in\Tau_{\Omega_i}$, and from \eqref{eq:localMassConservation}  we see that $U_{\mathrm{LC}}^\perp = \prod_{i=1}^m \mathbb{P}_0(\Tau_{\Omega_i})$. 

We will consider the slightly weaker case, where \eqref{eq:localMassConservation} is only required to hold for all ``non-Dirichlet boundary'' elements, that is for all elements where $\partial K\cap \partial_D\Omega=\emptyset$. This is sufficient for the results from Lemma~\ref{lemma:localPoincare} to be extendable to the grid level by considering the space $W_\Omega = \prod_{i=1}^m \prod_{K\in\Tau_{\Omega_i}} \mathring{H}^1(K)$, where $\mathring{H}^1(K)$ is defined in~\eqref{eq:zeroMean}.

Lemma \ref{lemma:localPoincare} now applies without modification, and weights ${\mu}(\vecx) \geq C^{-1}_{K}$ for $\vecx \in K$ are therefore in $\calC_W$. Moreover, thanks to convexity of simplicial grid elements, the local permeability-weighted Poincar{\'e}-Friedrichs constants are now fully computable. This allows us to bound $\eta_{\mathrm{R},\Omega}$ as follows: 
\begin{align}
    \eta_{\mathrm{R},\Omega}^2 &\leq \sum_{i=1}^m \sum_{K\in\Tau_{\Omega_i}} \frac{h_K^2}{\pi^2 c_{K}^2} \norm{f_i - \nabla_i\cdot\Bigg(
\vecv_{0,i} + \sum_{j\in\slow}\mathcal{R}_j\nu_{j}
\Bigg)
+ \sum_{j\in\shigh}\nu_{j}}_K^2 \nonumber \\
& = \sum_{\substack{i=1}}^m \sum_{K\in\Tau_{\Omega_i}} \eta_{\mathrm{R},K;\mathrm{LC}}^2 = \eta^2_{\mathrm{R},\Omega;\mathrm{LC}},\label{eq:eta_LC}
\end{align}
where $\eta_{\mathrm{R},K;\mathrm{LC}}$ are the local residual estimators for locally mass-conservative approximations. Using the above results, the majorant for locally mass-conservative approximations reads
\begin{equation}
    \mathcal{M}_\mathrm{LC}(\gothq,\gothv,\gothf) = \eta_{\mathrm{DF}}(\gothq,\gothv) + \eta_{\mathrm{R};\mathrm{LC}}(\gothv,\gothf).\label{eq:M_LC}
\end{equation}
The local residual estimates $\eta_{\mathrm{R},\Omega;\mathrm{LC}}$ correspond to the ones previously obtained by~\cite{vohralik2010unified, ernGuaranteed2017} for mono-dimensional problems subject to a flux equilibration step. Since $C_{K}\leq C_{\Omega_i}$, then, as before, whenever the assumptions of this section are satisfied, it holds that $\mathcal{M}\leq\mathcal{M}_\mathrm{LC}\leq\mathcal{M}_\mathrm{SC}\leq\mathcal{M}_\mathrm{NC}$.

\begin{remark}[Fully computable residual estimators] Unlike estimators obtained with residual methods (containing unknown constants~\cite{babuska1978error,kelly1983aposteriori}) or a purely functional approach such as in Sections \ref{sec:nocon} and \ref{sec:submasscon} (containing constants that are generally difficult to determine~\cite{repin2008posteriori}), estimators such as~\eqref{eq:eta_LC} contain only known local constants depending on the mesh size and material parameters. This justifies the claim that these estimators are \textit{fully computable}.
\end{remark}

\subsubsection{Exact mass-conservation\label{sec:EC}}

Methods with local mass conservation, as discussed in the previous section, when applied to problems where the RHS data $\gothf$ is zero or piecewise constant, can then often be verified to have an exact (pointwise) conservation property. We indicate this case by the abbreviation ``EC'', for which $\gothf = \mixdiv\gothv$, so that $U_{\mathrm{EC}} = 0$ and $U_{\mathrm{EC}}^{\perp} = L^2(\Omega)$. Now, $\pi_{W} = 0$ and $W = 0$. Thus, any finite weights $\mu$ are admissible, yet the choice is immaterial since the residual term $\norm{\mu^{-1}(\gothf-\mixdiv\gothv)}_{\Omega}$ always evaluates to zero. Consequently, only diffusive-type errors are present in the \textit{a posteriori} estimation, and the majorant takes the form
\begin{equation}
    \mathcal{M}_\mathrm{EC}(\gothq,\gothv) = \eta_{\mathrm{DF}}(\gothq,\gothv).\label{eq:M_EC}
\end{equation}
This case can also be seen as the limiting case of local mass conservation for a family of grid partitions where $h_K\rightarrow 0$.

\subsubsection{Summary of majorants and subdomain errors}

With the obtained majorants, we can define the corresponding upper bounds for the errors of the primal, dual, and primal-dual variables.

\begin{definition} Let $\alpha = \mathrm{NC, SC, LC, EC}$, corresponding to the flux conformity spaces $U_\alpha$ discussed in the preceding sections. Then, in view of the results from Theorem~\ref{thm:Abstract} and the majorants~\eqref{eq:M_NC}, \eqref{eq:M_SC}, \eqref{eq:M_LC}, and \eqref{eq:M_EC}, the upper bounds for the error in the primal, dual, and primal-dual pair, for arbitrary approximations $\gothq\in\mixHIzero+\gothg$ and $\gothv\in\mixHdivUa{\alpha}$, are
\begin{equation}
    \mathcal{M}^{\oplus}_{\gothp;\alpha} := \mathcal{M}_{\alpha}, \qquad \mathcal{M}^{\oplus}_{\gothu;\alpha} := \mathcal{M}_{\alpha}, \qquad \mathcal{M}^{\oplus}_{\gothp,\gothu;\alpha} := 2\mathcal{M}_{\alpha} + \eta_{\mathrm{R};\alpha}.
    \label{def:majorantBounds}
\end{equation}
while the lower bound for the error in the primal-dual pair is
\begin{equation}
    \mathcal{M}^{\ominus}_{\gothp,\gothu;\alpha} := \mathcal{M}_{\alpha}.
\end{equation}
\end{definition}

It is our interest not only to measure local errors, but also to distinguish between subdomain and interface errors. This motivates the definition of the following errors estimators.  

\begin{definition}[Subdomain and interface error indicators] Let $\alpha=\mathrm{EC, LC, SC, NC}$. Then, we will denote by $\varepsilon_{\Omega_i;\alpha}$ and $\varepsilon_{\Gamma_{j}}$ the subdomain and interface error indicators, defined by
\begin{align*}
    \varepsilon_{\Omega_i;\alpha}^2 &:= \varepsilon_{\mathrm{DF},\Omega_i}^2 + \varepsilon_{\mathrm{R},\Omega_i;\alpha}^2 := \sum_{K\in\Tau_{\Omega_i}} \eta^2_{\mathrm{DF}_{\parallel},K} + \sum_{K\in\Tau_{\Omega_i}} \eta^2_{\mathrm{R},K;\alpha},\\
    \varepsilon^2_{\Gamma_{j}} &:= \varepsilon_{\mathrm{DF},\Gamma_{j}}^2 := \sum_{K\in\Tau_{\Gamma_{j}}} \eta^2_{\mathrm{DF}_\perp,K}.
\end{align*}
\label{def:errors}
\end{definition}

We emphasize that while the majorants provide guaranteed bounds, the subdomain and interface error indicators can only be expected to correlate with the error.  

\section{Concrete bounds for locally mass-conservative approximations\label{sec:concreteBounds}}

In this section, we will make the evaluation of the bounds concrete by providing explicit approximations to~\eqref{eq:mD_weakDual} using the lowest-order mixed-finite element method (MFEM).

\subsection{Grid partitions\label{subsec:gridPartitions}}

Ultimately, \textit{a posteriori} estimates are primarily applied to approximations that are defined on computational grids. We therefore, in this section, summarize the relevant notation for grids and the mapping operators between subdomains and interfaces. 

Let us start by defining the partitions of the domains of interest. To this aim, denote by $\Tau_{\Omega_i}$, $\Tau_{\Gamma_j}$, and $\Tau_{\partial_i\Omega_j}$ the partitions of $\Omega_i$, $\Gamma_j$, and $\partial_j\Omega_i$, respectively. Moreover, let $\Tau_{\Omega} = \cup_{i=1}^m \Tau_{\Omega_i}$, $\Tau_{\Gamma} = \cup_{j=1}^M \Tau_{\Gamma_j}$, and $\Tau_{\partial_I\Omega} = \cup_{i=1}^m \cup_{j\in\slow} \partial_j\Omega_i$ represent the union of all subdomain, mortar, and internal boundary grids. 

Here, we only consider simplicial partitions. In particular, we require all elements $K\subset\Omega_i$ to be strictly non-overlapping simplices of dimension $d_K=d_i$. We use $h_{K}$ to denote the diameter of $K$, and define $h_{\Omega_i} = \max_{h_K} \Tau_{\Omega_i}$, $h_{\Gamma_j} = \max_{h_K} \Tau_{\Gamma_j}$, and $h_{\partial_j\Omega_i} = \max_{h_K} \Tau_{\partial_j\Omega_i}$.

\begin{figure}[tbp]
    \centering
    \includegraphics[width=0.8\textwidth]{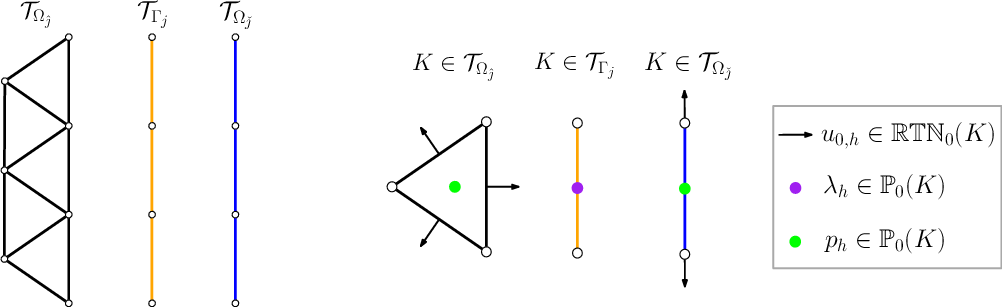}
    \caption{Left: Matching coupling between the grids $\Tau_{\Omega_\hatj}$, $\Tau_{\Gamma_j}$, and $\Tau_{\Omega_\checkj}$. Right: Degrees of freedom involved in the coupling between a 2d higher-dimensional cell, a 1d mortar-cell, and a 1d lower-dimensional cell. Locally, tangential fluxes are approximated using $\mathbb{RTN}_0(K)$, whereas mortar fluxes and pressures using $\mathbb{P}_0(K)$.
    \label{fig:matching_grids}}
\end{figure}

We will not at this point place any conditions on the grid partitions, although several aspects of this will be advantageous from the perspective of computation.

\subsection{Finite element spaces and the approximated problem\label{sec:fem_spaces}}

Let us introduce the finite element spaces necessary to write the approximated problem. We start by defining a local space for the approximated pressures, mortar fluxes, and tangential fluxes. They are given, respectively by 
\begin{alignat*}{2}
    Q_{h,i} &:= \left\lbrace q_{h,i} \in L^2(\Omega_i) : q_{h,i} \rvert_K \in \mathbb{P}_0 (K) \,\, \forall\, K \in \Tau_{\Omega_i} \right\rbrace, \quad&&d_i \in \{0, \ldots, n\}, \\
    \Lambda_{h,j} &:= \left\lbrace \nu_{h,j} \in L^2 (\Gamma_{j}) : \nu_{h,j}\rvert_K \in \mathbb{P}_0(K) \,\, \forall \, K \in \Tau_{\Gamma_{j}} \right\rbrace, \quad&&d_j \in \{0, \ldots, n-1\}, \\
    {V}_{h,i} &:= \left\lbrace \vecv_{h,i} \in H(\mathrm{div}; \Omega_i) : \vecv_{h,i} \rvert_K \in \mathbb{RTN}_0(K) \,\, \forall \, K \in \Tau_{\Omega_i} \right\rbrace, \quad&&d_i \in \{1, \ldots, n\},
\end{alignat*}
where $\mathbb{P}_0$ and $\mathbb{RTN}_0$ denote the spaces of constants and lowest-order Raviart-Thomas(-N{\'e}d{\'e}lec) spaces of vector functions ~\cite{raviart1977mixed, nedelec1980mixed}. See also Figure~\ref{fig:matching_grids} for the degrees of freedom involved in the generic coupling between a (higher-dimensional) triangle, a mortar line segment, and a (lower-dimensional) line segment.

The composite space for the approximated mD~pressure ${Q}_h \subset \mixLIIscalar$ and the approximated mD~flux ${X}_h \subset \mixHdiv$ are defined respectively by
\begin{equation}
    {Q}_h := \prod_{i=1}^m Q_{h,i} \qquad\mathrm{and}\qquad {X}_h := \prod_{i=1}^m \left( H_0(\mathrm{div}; \Omega_i)\cap {V}_{h,i}  \times \prod_{j \in \slow} \mathcal{R}_{h,j} \, \Lambda_{h,j}\right) .
\end{equation}
While not strictly necessary from a theoretical perspective, in the discrete setting, it is often useful to choose a finite-dimensional reconstruction operator based on the discrete spaces, and we allow for this through the notation $\mathcal{R}_{h,j} : \Lambda_{h,j}\rightarrow H(\mathrm{div}; \Omega_i)$, which in practice is often further restricted to $\mathcal{R}_{h,j} : \Lambda_{h,j}\rightarrow {V}_{h,\hatj}$. Such discrete reconstruction operators are natural for matching grids, and can also be constructed in the more general case of non-matching grids, see e.g.,~\cite{boon2018robust, boon2020nonmatching, arbogast2000nonmatching}.
Here $\Pi_h : \Lambda_{h,j}\rightarrow \tilde{\Lambda}_{h,j}$ is the $L^2$ projection from the mortar grid on $\Gamma_j$ to the boundary simplicial partition of $\Omega_\hatj$. 

We have now all the elements necessary to write the finite-dimensional approximation to the dual mixed problem \eqref{eq:mD_weakDual}. 

\begin{definition}[Approximated mD~dual mixed formulation] Find $(\gothu_h,\gothp_h) \in {X}_h\times{Q}_h$ such that
\begin{subequations}
    \begin{alignat}{2}
        \inner{\mixperm^{-1}\gothu_h}{\gothv_h}_{\Omega, \Gamma} - \inner{\gothp_h}{\mixdiv\gothv_h}_{\Omega} &= \innerangle{\gothg_D}{\mixTraceD{\gothv_h}}_{\partial_D\Omega} &&\qquad\forall\,\gothv_h\in{X}_h, \label{eq:mD_weakDualDiscrete1}\\
        \inner{\mixdiv\gothu_h}{\gothq_h}_{\Omega} &= \inner{\gothf}{\gothq_h}_{\Omega} &&\qquad\forall\,\gothq_h\in{Q}_h. \label{eq:mD_weakDualDiscrete2}
    \end{alignat}
\label{eq:mD_weakDualDiscrete}
\end{subequations}
\end{definition}

Due to the presence of the discrete reconstruction operator, this approximation is conforming whenever $\Lambda_{h,j}=\tilde{\Lambda}_{h,j}$, i.e., for matching grids. For non-matching grids, the approximation is still convergent, subject to normal conditions on the mortar grids~\cite{boon2018robust}.

\begin{remark}[Conservation properties]
\label{rem:cons_prop} Whenever equation \eqref{eq:mD_weakDualDiscrete2} is satisfied exactly, then equation \eqref{eq:localMassConservation} holds, and we have local mass conservation for matching grids. Thus, the fluxes lie in the smaller space ${X}_h\cap\mixHdivX{Q_{h,i}^\perp}$, and the results from section \ref{sec:LMC} apply. Furthermore, if $f_i\in Q_{h,i}$ and $\mathcal{R}_{h,j} : \Lambda_{h,j}\rightarrow {V}_{h,\hatj}$, then the projection of the source term, and hence the residual error, onto $Q_{h,i}^\perp$ vanishes. Thus, the local conservation is verified to be pointwise, the fluxes lie in ${X}_h\cap\mixHdivX{0}$ and the results from Section \ref{sec:EC} apply.
\end{remark}

\begin{remark}[Well-posedness and a priori estimates] The stability and \textit{a priori} approximation properties of the finite-dimensional system given in \eqref{eq:mD_weakDualDiscrete} has been previously established~\cite{boon2018robust}.
\end{remark}

\subsection{Pressure reconstruction}
Recall that Theorem~\ref{thm:Abstract} requires any approximation to the mD~flux to be in $\mixHdiv$, whereas approximations to the mD~pressure must lie in $\mixHIzero + \gothg$. By the condition that $\gothu_h \in {X}_h\subset\mixHdiv$, the solution of equations \eqref{eq:mD_weakDualDiscrete} by definition satisfy the first condition.
On the other hand, the approximated mD~pressure $\gothp_h$ is only in $\mixLIIscalar$. We therefore need to enhance the regularity of the approximated pressure and thus obtain a reconstructed pressure.

\begin{definition}[Reconstructed pressure\label{def:pressureReconstruction}] We will call reconstructed pressure $\gothpp$ to any function constructed from the mD~pair $(\gothp_h,\gothu_h) \in \mixLIIscalar\times\mixHdiv$ such that
\begin{equation}
    \gothpp \in \mixHIzero+\gothg.
\end{equation}
\end{definition}

\begin{remark}[On potential reconstruction]
\label{rem:reconstruction} Several techniques for obtaining $\gothpp$ are available in the literature. Arguably, the simplest option is to perform an average of the $\mathbb{P}_0(K)$ pressures on local patches and from there construct local affine $\mathbb{P}_1(K)$ functions~\cite{repin2009aposterior}. Other techniques aim at solving first a local Neumann problem to obtain a $\mathbb{P}_2(K)$ post-processed pressure, and then apply interpolation techniques to get energy-conforming potentials~\cite{ainsworth2007mixed, ern2010posteriori, vohralik2010unified, ahmed2019adaptive, ahmed2020adaptive}. Any of these choices are compatible with the bounds derived herein.
\end{remark}

\begin{remark}[Computable estimates] Computable versions of the majorants are now readily available by setting $(\gothq,\gothv) = (\gothpp,\gothu_h)$ in \eqref{eq:M_NC}, \eqref{eq:M_SC}, \eqref{eq:M_LC}, and \eqref{eq:M_EC}.
\end{remark}

\begin{remark}[Other locally mass-conservative methods]
\label{rem:loc_mass}
In addition to the MFEM scheme of the lowest-order (RT0-P0), other flux-based numerical methods such as the Mixed Virtual Element Method (MVEM) \cite{daVeiga2016MVEM, fumagalli2019dual} and Cell Centered Finite Volume Methods (CCFVM),  including the Two-Point Flux Approximation (TPFA) and the Multi-Point Flux Approximation (MPFA) \cite{aavatsmark2002introduction, nordbotten2021introduction}, can be analyzed with our framework provided that the fluxes are interpolated in ${X}_{h}$ and the pressures reconstructed as indicated above. For methods without an explicit flux representation, an additional flux reconstruction step may be needed. 
\end{remark}

\begin{remark}[Superconvergence of the residual estimators\label{rem:superconvergence}] Due to Remark \ref{rem:cons_prop}, the residual estimators $\eta_{\mathrm{R},K,\mathrm{LC}}$ are superconvergent for lowest-order locally mass-conservative approximations. This property is guaranteed since: (1) local Poincar{\'e} constants decay as $\mathcal{O}(h_K)$ for simplicial elements and (2) the norm of the residual $\norm{f_i - \nabla_{i} \cdot \vecv_i + \textstyle{\sum_{j\in\hat{S}_i}{\nu_j}}}_K$ also decays as $\mathcal{O}(h_K)$~\cite{boffi2013mixed}; leading to an overall rate of $\mathcal{O}(h_K^2)$. 
\end{remark}

\section{Numerical validations\label{sec:num_val}}

\begin{figure}
	\centering
	\includegraphics[height=0.225\textheight]{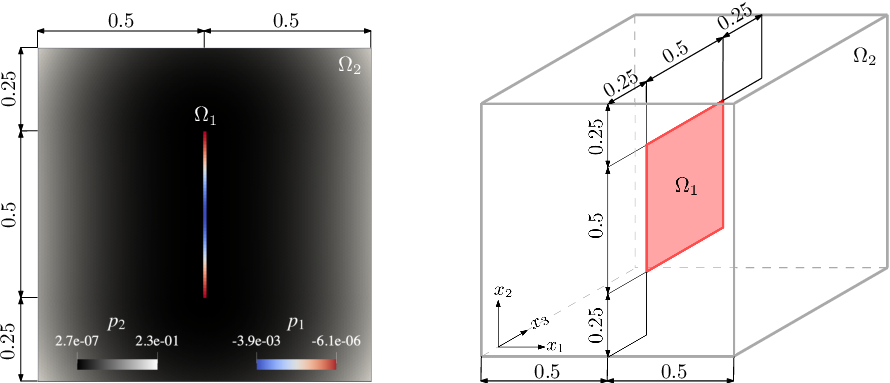}
	\caption{Geometric setups used for the numerical validations. Left: A $1d$ fracture embedded in a $2d$ matrix and the exact pressure solution. Right: A $2d$ fracture embedded in a $3d$ matrix. \label{fig:validation_domains}}
\end{figure}

In this section, we test the performance of our estimators by conducting an efficiency analysis using four different numerical methods, namely those mentioned in Remark \ref{rem:loc_mass}: RT0-P0, MVEM-P0, MPFA, and TPFA. The numerical examples are implemented in the Python-based open-source software \texttt{PorePy}~\cite{keilegavlen2020porepy}, using the extension package \texttt{mdestimates}~\cite{mdestimates}, which includes the scripts of all numerical examples considered here. In these numerical validations, we only consider matching grids, and use a low-order pressure reconstruction (recall Remark \ref{rem:reconstruction} for further discussion).  

We validate the \textit{a posteriori} bounds and assess their efficiency on a 1d/2d problem  (Section~\ref{sec:2d-validation}) and a 2d/3d problem (Section~\ref{sec:3d-validation}), both with manufactured solutions. The geometric configuration for both problems is shown in Figure~\ref{fig:validation_domains}. Let us denote the fracture as $\Omega_1$, the matrix as $\Omega_2$, the left interface as $\Gamma_{1}$, and the right interface as $\Gamma_{2}$. Further, assume the existence of an exact, smooth pressure $p_2(\vecx)$ in $\Omega_2$. Refer to Table~\ref{tab:exactSolution2D} and Table~\ref{tab:exactSolution3D} from the Appendix~\ref{sec:exactSolutions} for the analytical expressions of all variables of interest.

\subsection{Efficiency indices\label{sec:efficiencyIndices}}

Efficiency indices are used to assess the performance of the approximations when exact solutions are available. They are defined as the ratio between the estimated and the exact errors. Here, we consider the following efficiency indices.

\begin{definition}[Efficiency indices\label{def:indices}] Let $\alpha = \mathrm{NC, SC, LC, EC}$ and let $\gothp \in \mixHIzero + \gothg$ and $\gothu \in \mixHdiv$ be the solutions to~\eqref{eq:mD_weakPrimal} and~\eqref{eq:mD_weakDual}, respectively. Then, in view of Theorem~\ref{thm:Abstract}, the efficiency indices for the primal, dual, and primal-dual pair, for arbitrary approximations $\gothq\in\mixHIzero+\gothg$ and $\gothv\in\mixHdivUa{\alpha}$, are
\begin{equation}
    I_{\gothp;\alpha}(\gothq) := \frac{\mathcal{M}^{\oplus}_{\gothp;\alpha}}{\tnorm{\gothp-\gothq}}, \qquad  I_{\gothu;\alpha}(\gothv) := \frac{\mathcal{M}^{\oplus}_{\gothu;\alpha}}{\tnormstar{\gothu-\gothv}}, \qquad I_{\gothp,\gothu;\alpha}(\gothq,\gothv) := \frac{\mathcal{M}^{\oplus}_{\gothp,\gothu;\alpha}}{\fullnorm{\gothp-\gothq, \gothu-\gothv}}.
\end{equation}
\end{definition}

\begin{remark}\label{rem:eff_ind}
Optimal efficiency indices (equal to 1) are obtained when the approximations match the exact solutions. Moreover, in general the efficiency indices satisfy the bounds: 
\begin{equation}
    1\leq I_{\gothp;\alpha}(\gothq), \qquad 1\leq I_{\gothu;\alpha}(\gothv), \qquad 1\leq I_{\gothp,\gothu;\alpha}(\gothq,\gothv)\leq  \frac{\mathcal{M}^{\oplus}_{\gothp,\gothu;\alpha}}{\mathcal{M}^{\ominus}_{\gothp,\gothu;\alpha}} = 2 +  \frac{\eta_{\mathrm{R};\alpha}}{\mathcal{M}_\alpha}
\end{equation}
For the final term, we note that since $\eta_{\mathrm{R};\alpha}\leq\mathcal{M}_\alpha $, then for $\alpha = \mathrm{NC, SC}$ the total efficiency index satisfies $I_{\gothp,\gothu;\alpha}\leq 3$, while for local conservation $I_{\gothp,\gothu;\mathrm{LC}}\leq 2 + O(h^2)$ and finally for exact conservation $I_{\gothp,\gothu;\mathrm{E C}}\leq 2$.
\end{remark}

\subsection{Two-dimensional validation\label{sec:2d-validation}}

\begin{table}
\caption{Two-dimensional validation: Majorants and efficiency indices.\label{tab:2d-validation_b}}
\centering
\begin{adjustbox}{width=\textwidth}
    \begin{tabular}{r r a b a b a b a b a b}
    \toprule
    & $h_{\mathrm{coup}}$ & $\mathcal{M}^{\oplus}_{\gothp;\mathrm{NC}}$ &  $\mathcal{M}^{\oplus}_{\gothp;\mathrm{LC}}$ & $\mathcal{M}^{\oplus}_{\gothp,\gothu;\mathrm{NC}}$ & $\mathcal{M}^{\oplus}_{\gothp,\gothu;\mathrm{LC}}$ & $I_{\gothp;\mathrm{NC}}$ & $I_{\gothp;\mathrm{LC}}$ &  $I_{\gothu;\mathrm{NC}}$ & $I_{\gothu;\mathrm{LC}}$ & $I_{\gothp,\gothu;\mathrm{NC}}$ & $I_{\gothp,\gothu;\mathrm{LC}}$ \\
    \midrule
    \multirow{4}{*}{\rotatebox[origin=c]{90}{RT0-P0}} & 0.05 & 5.86e-02 & 4.36e-02 & 1.33e-01 & 8.83e-02 & 1.46 & 1.08 & 4.09 & 3.04 & 1.89 & 1.59 \\
    & 0.025 & 3.01e-02 & 2.17e-02 & 6.89e-02 & 4.38e-02 & 1.49 & 1.07 & 4.18 & 3.02 & 1.91 & 1.58  \\
    & 0.0125 & 1.52e-02 & 1.08e-02 & 3.48e-02 & 2.17e-02 & 1.50 & 1.07 & 4.22 & 3.00 & 1.92 & 1.57 \\
    & 0.00625 & 7.65e-03 & 5.37e-03 & 1.76e-02 & 1.08e-02 & 1.52 & 1.07 & 4.25 & 2.98 & 1.93 & 1.57 \\
    \midrule
    \multirow{4}{*}{\rotatebox[origin=c]{90}{MVEM-P0}} & 0.05 & 6.18e-02 & 4.68e-02 & 1.40e-01 & 9.47e-02 & 1.42 & 1.07 & 4.31 & 3.26 & 1.89 & 1.60 \\ 
    & 0.025 & 3.10e-02 & 2.27e-02 & 7.08e-02 & 4.56e-02 & 1.46 & 1.07 & 4.31 & 3.15 & 1.91 & 1.59 \\
    & 0.0125 & 1.54e-02 & 1.10e-02 & 3.53e-02 & 2.22e-02 & 1.49 & 1.07 & 4.29 & 3.07 & 1.92 & 1.58 \\
    & 0.00625 & 7.72e-03 & 5.44e-03 & 1.77e-02 & 1.09e-02 & 1.51 & 1.06 & 4.28 & 3.02 & 1.92 & 1.57 \\
    \midrule
    \multirow{4}{*}{\rotatebox[origin=c]{90}{MPFA}} & 0.05 & 5.91e-02 & 4.41e-02 & 1.34e-01 & 8.93e-02 & 1.46 & 1.09 & 4.12 & 3.07 & 1.89 & 1.59 \\
    & 0.025 & 3.03e-02 & 2.19e-02 & 6.92e-02 & 4.41e-02 & 1.49 & 1.08 & 4.20 & 3.04 & 1.91 & 1.58 \\
    & 0.0125 & 1.52e-02 & 1.08e-02 & 3.49e-02 & 2.18e-02 & 1.50 & 1.07 & 4.23 & 3.01 & 1.92 & 1.57 \\
    & 0.00625 & 7.66e-03 & 5.38e-03 & 1.76e-02 & 1.08e-02 & 1.52 & 1.07 & 4.25 & 2.99 & 1.93 & 1.57 \\
    \midrule
    \multirow{4}{*}{\rotatebox[origin=c]{90}{TPFA}} & 0.05 & 6.67e-02 & 5.17e-02 & 1.50e-01 & 1.04e-01 & 1.54 & 1.19 & 3.09 & 2.39 & 1.84 & 1.58 \\
    & 0.025 & 3.74e-02 & 2.90e-02 & 8.35e-02 & 5.83e-02 & 1.68 & 1.31 & 2.36 & 1.83 & 1.78 & 1.52 \\
    & 0.0125 & 2.64e-02 & 2.20e-02 & 5.73e-02 & 4.41e-02 & 1.82 & 1.52 & 1.64 & 1.36 & 1.63 & 1.44 \\
    & 0.00625 & 1.37e-02 & 1.15e-02 & 2.98e-02 & 2.30e-02 & 1.64 & 1.37 & 1.82 & 1.52 & 1.64 & 1.44 \\
    \bottomrule
    \end{tabular}
\end{adjustbox}
    \vspace{1ex}
    
    {\scriptsize\raggedright The results for $\mathcal{M}^{\oplus}_{\gothu;\mathrm{NC}}$ and $\mathcal{M}^{\oplus}_{\gothu;\mathrm{LC}}$ are omitted since they are equal to $\mathcal{M}^{\oplus}_{\gothp;\mathrm{NC}}$ and $\mathcal{M}^{\oplus}_{\gothp;\mathrm{LC}}$. \par}
\end{table}

\begin{table}
    \caption{Two-dimensional validation: Subdomain and interface errors.\label{tab:2d-validation_a}}
    \begin{adjustbox}{width=\textwidth}
    \begin{tabular}{r r c a b c a b c c c}
    \toprule
    & $h_{\mathrm{coup}}$ & $\varepsilon_{\mathrm{DF},\Omega_2}$ & $\varepsilon_{\mathrm{R},\Omega_2;\mathrm{NC}}$  & $\varepsilon_{\mathrm{R},\Omega_2;\mathrm{LC}}$ & $\varepsilon_{\mathrm{DF},\Omega_1}$ & $\varepsilon_{\mathrm{R},\Omega_1;\mathrm{NC}}$ & $\varepsilon_{\mathrm{R},\Omega_1;\mathrm{LC}}$ & $\varepsilon_{\mathrm{DF},\Gamma_1}$ & $\varepsilon_{\mathrm{DF},\Gamma_2}$ \\
    \midrule
    \multirow{4}{*}{\rotatebox[origin=c]{90}{RT0-P0}} & 0.05 & 4.24e-02 & 1.41e-02 & 1.00e-03 & 2.26e-03 & 7.99e-03 & 5.65e-04 & 1.89e-04 & 1.89e-04 \\
    & 0.025 & 2.14e-02 & 7.73e-03 & 3.02e-04 & 1.14e-03 & 4.01e-03 & 1.42e-04 & 9.03e-05 & 9.15e-05  \\
    & 0.0125 & 1.07e-02 & 4.00e-03 & 7.28e-05 & 5.70e-04 & 2.01e-03 & 3.55e-05 & 4.41e-05 & 4.41e-05 \\
    & 0.00625 & 5.34e-03 & 2.07e-03 & 1.91e-05 & 2.85e-04 & 1.00e-03 & 8.87e-06 & 2.20e-05 & 2.20e-05 \\
    \midrule
    \multirow{4}{*}{\rotatebox[origin=c]{90}{MVEM-P0}} & 0.05 & 4.55e-02 & 1.41e-02 & 1.00e-03 & 3.25e-03 & 7.99e-03 & 5.65e-04 & 2.52e-04 & 2.52e-04 \\
    & 0.025 & 2.23e-02 & 7.73e-03 & 3.02e-04 & 1.32e-03 & 4.01e-03 & 1.42e-04 & 1.00e-04 & 1.03e-04 \\
    & 0.0125 & 1.09e-02 & 4.00e-03 & 7.28e-05 & 5.98e-04 & 2.01e-03 & 3.55e-05 & 4.50e-05 & 4.50e-05 \\
    & 0.00625 & 5.41e-03 & 2.07e-03 & 1.91e-05 & 2.89e-04 & 1.00e-03 & 8.87e-06 & 2.21e-05 & 2.21e-05 \\
    \midrule
    \multirow{4}{*}{\rotatebox[origin=c]{90}{MPFA}} & 0.05 & 4.29e-02 & 1.41e-02 & 1.00e-03 & 2.52e-03 & 7.99e-03 & 5.65e-04 & 2.05e-04 & 2.05e-04 \\
    & 0.025 & 2.15e-02 & 7.73e-03 & 3.02e-04 & 1.18e-03 & 4.01e-03 & 1.42e-04 & 9.24e-05 & 9.35e-05 \\
    & 0.0125 & 1.07e-02 & 4.00e-03 & 7.28e-05 & 5.77e-04 & 2.01e-03 & 3.55e-05 & 4.44e-05 & 4.44e-05 \\
    & 0.00625 & 5.36e-03 & 2.07e-03 & 1.91e-05 & 2.86e-04 & 1.00e-03 & 8.87e-06 & 2.20e-05 & 2.20e-05 \\
    \midrule
    \multirow{4}{*}{\rotatebox[origin=c]{90}{TPFA}} & 0.05 & 5.04e-02 & 1.41e-02 & 1.00e-03 &2.52e-03 & 7.99e-03 & 5.65e-04 & 1.87e-04 & 1.89e-04 \\
    & 0.025 & 2.86e-02 & 7.73e-03 & 3.02e-04 & 1.18e-03 & 4.01e-03 & 1.42e-04 & 9.42e-05 & 9.23e-05 \\
    & 0.0125 & 2.19e-02 & 4.00e-03 & 7.28e-05 & 5.77e-04 & 2.01e-03 & 3.55e-05 & 4.47e-05 & 4.46e-05 \\
    & 0.00625 & 1.14e-02 & 2.07e-03 & 1.91e-05 & 2.86e-04 & 1.00e-03 & 8.87e-06 & 2.20e-05 & 2.21e-05 \\
    \bottomrule
    \end{tabular}
    \end{adjustbox}
\end{table}

For our first validation, we consider the 1d/2d case as shown in the left Figure~\ref{fig:validation_domains}. This validation has two purposes: 
(1) compare the majorants and efficiency indices obtained using global (no mass-conservation) and local (local mass-conservation) Poincar{\'e}-Friedrichs constants, and 
(2) show the different errors associated with subdomains and interfaces.

To this aim, we consider four levels of successively refined combinations of mesh sizes, characterized by  $h_{\mathrm{coup}} = h_{\partial_1\Omega_2} = h_{\Gamma_1} = h_{\Omega_1} = h_{\Gamma_{2}} =  h_{\partial_2\Omega_2}$. The global Poincar{\'e} constant is obtained numerically by solving the associated eigenvalue problem (see e.g.,~\cite{paulyPoincare2020}), giving a value of $C_{\Omega, \Gamma} \approx 0.2251$. 

Majorants for the primal, dual, and primal-dual variables are shown in Table~\ref{tab:2d-validation_b}. We can see that all majorants reflect the convergence tendency of the numerical methods, and in particular (as is well-known), we identify that the TPFA approximation performs relatively poorly on this problem. As expected, the majorants obtained exploiting the local conservation properties of the methods are sharper than the ones obtained using global weights, both in absolute value and in terms of efficiency index.

Further inspection shows that efficiency indices lie within the expected bounds discussed in Remark \ref{rem:eff_ind}. In particular, efficiency indices for the primal variable using local weights are very accurate, and only a $\sim{}7\%$ deviation with respect to the actual error (for the finest grid) is observed in the case of RT0-P0, MVEM-P0, and MPFA. For TPFA, the efficiency index is worse, as a consequence of the flux approximation being worse. Efficiency indices for the dual variable are in general larger than the ones obtained for the primal variable; this is to be expected for mixed-dual approximations with the relatively simple pressure reconstruction, where the approximated fluxes have relatively good accuracy as compared to the reconstructed pressures. Finally, efficiency indices for the primal-dual variable are less than $2$ for all methods in consideration.

Considering now the local error indicators, shown in Table~\ref{tab:2d-validation_a}, we note that diffusive errors decrease linearly for the matrix, fracture, and interfaces. Likewise, residual errors for the matrix and fracture decrease linearly when the global Poincar{\'e}-Friedrichs constant is used. When the local Poincar{\'e}-Freidrich constants are used, the residual estimators for the matrix and the fracture decrease quadratically, which goes in agreement with the super-convergent properties discussed in Remark~\ref{rem:superconvergence}.

\subsection{Three-dimensional validation\label{sec:3d-validation}}

\begin{table}
    \caption{Three-dimensional validation: Majorants and efficiency indices.\label{tab:3d-validation_b}}
    \begin{adjustbox}{width=\textwidth}
    \begin{tabular}{r r a b a b a b a b a b}
    \toprule
    & $h_{\mathrm{coup}}$ & $\mathcal{M}^{\oplus}_{\gothp;\mathrm{NC}}$ &  $\mathcal{M}^{\oplus}_{\gothp;\mathrm{LC}}$ & $\mathcal{M}^{\oplus}_{\gothp,\gothu;\mathrm{NC}}$ & $\mathcal{M}^{\oplus}_{\gothp,\gothu;\mathrm{LC}}$ & $I_{\gothp;\mathrm{NC}}$ & $I_{\gothp;\mathrm{LC}}$ &  $I_{\gothu;\mathrm{NC}}$ & $I_{\gothu;\mathrm{LC}}$ & $I_{\gothp,\gothu;\mathrm{NC}}$ & $I_{\gothp,\gothu;\mathrm{LC}}$ \\
    \midrule
    \multirow{4}{*}{\rotatebox[origin=c]{90}{RT0-P0}} & 0.2625 & 2.85e-01 & 2.36e-01 & 6.21e-01 & 4.73e-01 & 1.25 & 1.03 & 4.14 & 3.43 & 1.72 & 1.43 \\
    & 0.1720 & 1.94e-01 & 1.62e-01 & 4.20e-01 & 3.24e-01 & 1.23 & 1.03 & 4.22 & 3.53 & 1.73 & 1.50  \\
    & 0.0827 & 1.07e-01 & 8.69e-02 & 2.33e-01 & 1.74e-01 & 1.28 & 1.04 & 3.61 & 2.94 & 1.70 & 1.48 \\
    & 0.0418 & 5.62e-02 & 4.58e-02 & 1.23e-01 & 9.16e-02 & 1.25 & 1.02 & 3.16 & 2.58 & 1.63 & 1.43 \\
    \midrule
    \multirow{4}{*}{\rotatebox[origin=c]{90}{MVEM-P0}} & 0.2625 & 2.89e-01 & 2.40e-01 & 6.28e-01 & 4.80e-01 & 1.24 & 1.03 & 4.19 & 3.48 & 1.72 & 1.44 \\ 
    & 0.1720 & 1.96e-01 & 1.64e-01 & 4.24e-01 & 3.28e-01 & 1.23 & 1.03 & 4.26 & 3.57 & 1.73 & 1.50  \\
    & 0.0827 & 1.08e-01 & 8.80e-02 & 2.35e-01 & 1.76e-01 & 1.27 & 1.04 & 3.65 & 2.98 & 1.70 & 1.48 \\
    & 0.0418 & 5.66e-02 & 4.62e-02 & 1.24e-01 & 9.23e-02 & 1.25 & 1.02 & 3.18 & 2.60 & 1.63 & 1.44 \\
    \midrule
    \multirow{4}{*}{\rotatebox[origin=c]{90}{MPFA}} & 0.2625 & 2.90e-01 & 2.40e-01 & 6.29e-01 & 4.82e-01 & 1.25 & 1.03 & 4.08 & 3.38 & 1.72 & 1.43 \\
    & 0.1720 & 1.98e-01 & 1.66e-01 & 4.28e-01 & 3.32e-01 & 1.23 & 1.03 & 4.22 & 3.54 & 1.73 & 1.50 \\
    & 0.0827 & 1.09e-01 & 8.90e-02 & 2.37e-01 & 1.78e-01 & 1.27 & 1.04 & 3.64 & 2.98 & 1.70 & 1.49 \\
    & 0.0418 & 5.69e-02 & 4.65e-02 & 1.24e-01 & 9.30e-02 & 1.25 & 1.02 & 3.18 & 2.60 & 1.63 & 1.44 \\
    \midrule
    \multirow{4}{*}{\rotatebox[origin=c]{90}{TPFA}} & 0.2625 & 3.84e-01 & 3.35e-01 & 8.17e-01 & 6.70e-01 & 1.24 & 1.08 & 2.13 & 1.86 & 1.48 & 1.28 \\
    & 0.1720 & 2.95e-01 & 2.63e-01 & 6.23e-01 & 5.27e-01 & 1.38 & 1.23 & 1.66 & 1.48 & 1.44 & 1.30 \\
    & 0.0827 & 2.22e-01 & 2.02e-01 & 4.63e-01 & 4.04e-01 & 1.62 & 1.48 & 1.40 & 1.28 & 1.45 & 1.35 \\
    & 0.0418 & 2.08e-01 & 1.97e-01 & 4.26e-01 & 3.95e-01 & 1.76 & 1.67 & 1.29 & 1.23 & 1.46 & 1.41 \\
    \bottomrule
    \end{tabular}
    \end{adjustbox}
    
    \vspace{1ex}
    {\scriptsize\raggedright The results for $\mathcal{M}^{\oplus}_{\gothu;\mathrm{NC}}$ and $\mathcal{M}^{\oplus}_{\gothu;\mathrm{LC}}$ are omitted since they are equal to $\mathcal{M}^{\oplus}_{\gothp;\mathrm{NC}}$ and $\mathcal{M}^{\oplus}_{\gothp;\mathrm{LC}}$. \par}
\end{table}

\begin{table}
    \caption{Three-dimensional validation: Subdomain and interface errors.\label{tab:3d-validation_a}}
    \begin{adjustbox}{width=\textwidth}
    \begin{tabular}{r r c a b c a b c c}
    \toprule
    & $h_{\mathrm{coup}}$ & $\varepsilon_{\mathrm{DF},\Omega_2}$ & $\varepsilon_{\mathrm{R},\Omega_2;\mathrm{NC}}$  & $\varepsilon_{\mathrm{R},\Omega_2;\mathrm{LC}}$ & $\varepsilon_{\mathrm{DF},\Omega_1}$ & $\varepsilon_{\mathrm{R},\Omega_1;\mathrm{NC}}$ & $\varepsilon_{\mathrm{R},\Omega_1;\mathrm{LC}}$ & $\varepsilon_{\mathrm{DF},\Gamma_1}$ & $\varepsilon_{\mathrm{DF},\Gamma_2}$ \\
    \midrule
    \multirow{4}{*}{\rotatebox[origin=c]{90}{RT0-P0}} & 0.2625 & 2.35e-01 & 4.73e-02 & 2.55e-02 & 4.67e-03 & 1.56e-02 & 6.70e-03 & 5.05e-03 & 5.00e-03  \\
    & 0.1720 & 1.62e-01 & 3.08e-02 & 1.10e-02 & 4.53e-03 & 9.01e-03 & 2.04e-03 & 1.37e-03 & 1.37e-03  \\
    & 0.0827 & 8.69e-02 & 1.91e-02 & 3.91e-03 & 2.72e-03 & 5.17e-03 & 6.70e-04 & 4.07e-04 & 4.09e-04 \\
    & 0.0418 & 4.58e-02 & 1.01e-02 & 1.06e-03 & 1.40e-03 & 2.64e-03 & 1.70e-04 & 1.08e-04 & 1.09e-04 \\
    \midrule
    \multirow{4}{*}{\rotatebox[origin=c]{90}{MVEM-P0}} & 0.2625 & 2.39e-01 & 4.73e-02 & 2.55e-02 & 5.78e-03 & 1.56e-02 & 6.70e-03 & 5.54e-03 & 5.49e-03 \\
    & 0.1720 & 1.64e-01 & 3.08e-02 & 1.10e-02 & 5.56e-03 & 9.01e-03 & 2.04e-03 & 1.50e-03 & 1.50e-03 \\
    & 0.0827 & 8.79e-02 & 1.91e-02 & 3.91e-03 & 3.05e-03 & 5.17e-03 & 6.70e-04 & 4.40e-04 & 4.41e-04 \\
    & 0.0418 & 4.61e-02 & 1.01e-02 & 1.06e-03 & 1.46e-03 & 2.64e-03 & 1.70e-04 & 1.13e-04 & 1.14e-04 \\
    \midrule
    \multirow{4}{*}{\rotatebox[origin=c]{90}{MPFA}} & 0.2625 & 2.40e-01 & 4.73e-02 & 2.55e-02 & 5.04e-03 & 1.56e-02 & 6.70e-03 & 5.93e-03 & 5.86e-03 \\
    & 0.1720 & 1.66e-01 & 3.08e-02 & 1.10e-02 & 4.86e-03 & 9.01e-03 & 2.04e-03 & 1.56e-03 & 1.55e-03 \\
    & 0.0827 & 8.89e-02 & 1.91e-02 & 3.91e-03 & 2.82e-03 & 5.17e-03 & 6.70e-04 & 4.61e-04 & 4.62e-04 \\
    & 0.0418 & 4.65e-02 & 1.01e-02 & 1.06e-03 & 1.41e-03 & 2.64e-03 & 1.70e-04 & 1.14e-04 & 1.16e-04 \\
    \midrule
    \multirow{4}{*}{\rotatebox[origin=c]{90}{TPFA}} & 0.2625 & 3.34e-01 & 4.73e-02 & 2.55e-02 & 4.88e-03 & 1.56e-02 & 6.70e-03 & 6.04e-03 & 5.11e-03 \\
    & 0.1720 & 2.63e-01 & 3.08e-02 & 1.10e-02 & 4.86e-03 & 9.01e-03 & 2.04e-03 & 1.35e-03 & 1.29e-03 \\
    & 0.0827 & 2.02e-01 & 1.91e-02 & 3.91e-03 & 2.85e-03  &5.17e-03 & 6.70e-04 & 4.50e-04 & 4.39e-04 \\
    & 0.0418 & 1.97e-01 & 1.01e-02 & 1.06e-03 & 1.46e-03 & 2.64e-03 & 1.70e-04 & 1.02e-04 & 1.02e-04 \\
    \bottomrule
    \end{tabular}
    \end{adjustbox}
\end{table}

For our next numerical validation, we employ the 2d/3d configuration from the right Figure~\ref{fig:validation_domains}. We repeat the same analysis from the previous section, and investigate four refinement levels. The mixed-dimensional Poincar{\'e} constant for this configuration corresponds to a value of $C_{\Omega,\Gamma} \approx 0.1838$. The results are shown in Table~\ref{tab:3d-validation_b} and
Table~\ref{tab:3d-validation_a}.
As in the previous validation, we can see that the majorants capture the local and global convergence tendency of all numerical methods. Again, RT0-P0, MVEM-P0, and MPFA give quite similar results, whereas TPFA showcase larger errors. As expected, efficiency indices again lie within the stated bounds from Remark \ref{rem:eff_ind}.

\section{Numerical applications\label{sec:numex_app}}

In this section, we apply our estimators to numerical approximations of challenging problems solving the equations of incompressible flow in fractured porous media. Importantly, since source terms are zero in both applications, by applying matching grids the residual errors are zero, and we are in the setting of having an exact conservation property from the numerical approximation. From Remark \ref{rem:eff_ind}, we then know that the efficiency index for the primal-dual error will be less than 2; even if the exact solution and error are both unknown.

\subsection{Two-dimensional application\label{sec:2dbench}}

\begin{table}
\caption{Error estimates for the two-dimensional application.\label{tab:bench2d}}
    \begin{adjustbox}{width=\textwidth}
    \begin{tabular}{r r c c c c c c c c}
    \toprule
    & Mesh & ${\varepsilon}_{\Omega^2;\mathrm{EC}}$ & $ {\varepsilon}_{\Omega^1;\mathrm{EC},\mathrm{C}}$ & ${\varepsilon}_{\Omega^1;\mathrm{EC},\mathrm{B}}$ & ${\varepsilon}_{\Gamma^1,\mathrm{C}}$ & ${\varepsilon}_{\Gamma^1,\mathrm{B}}$ & ${\varepsilon}_{\Gamma^0}$ & $\mathcal{M}^{\oplus}_{\gothp;\mathrm{EC}}$ & $\mathcal{M}^{\oplus}_{\gothp,\gothu;\mathrm{EC}}$ \\
    \midrule
    \multirow{3}{*}{\rotatebox[origin=c]{90}{RT0-P0~~}} & Coarse & 7.39e-01 & 2.93e-01 & 2.98e-04 & 3.13e+03 & 1.52e-01 & 2.24e+01 & 9.94e+02 & 1.99e+03 \\ [2mm]
    & Intermediate & 5.95e-01 & 1.90e-01 & 2.77e-04 & 1.95e+03 & 1.00e-01 & 1.79e+01 & 6.20e+02 & 1.24e+03 \\ [2mm]
    & Fine & 4.30e-01 & 1.07e-01 & 2.78e-04 & 9.79e+02 & 5.15e-02 & 1.22e+01 & 3.15e+02 & 6.30e+02 \\ 
    \midrule
    \multirow{3}{*}{\rotatebox[origin=c]{90}{MVEM-P0}} & Coarse & 7.29e-01 & 3.51e-01 & 1.44e-04 & 3.10e+03 & 1.46e-01 & 4.41e+01 & 9.84e+02 & 1.97e+03 \\ [2mm]
    & Intermediate & 5.91e-01 & 2.23e-01 & 1.27e-04 & 1.94e+03 & 9.43e-02 & 3.14e+01 & 6.17e+02 & 1.23e+03  \\ [2mm]
    & Fine & 4.28e-01 & 1.24e-01 & 1.18e-04 & 9.78e+02 & 4.80e-02 & 2.02e+01 & 3.15e+02 & 6.29e+02\\
    \midrule
    \multirow{3}{*}{\rotatebox[origin=c]{90}{MPFA~~}} & Coarse & 7.39e-01 & 3.13e-01 & 1.72e-04 & 3.03e+03 & 1.43e-01 & 3.39e+01 & 9.63e+02 & 1.93e+03 \\ [2mm]
    & Intermediate & 5.98e-01 & 2.01e-01 & 1.54e-04 & 1.89e+03 & 9.18e-02 & 2.55e+01 & 6.00e+02 & 1.20e+03 \\ [2mm]
    & Fine & 4.33e-01 & 1.12e-01 & 1.46e-04 & 9.49e+02 & 4.71e-02 & 1.68e+01 & 3.05e+02 & 6.10e+02 \\
    \midrule
    \multirow{3}{*}{\rotatebox[origin=c]{90}{TPFA~~}} & Coarse & 7.52e-01 & 3.05e-01 & 1.76e-04 & 3.19e+03 & 1.48e-01 & 3.67e+01 & 1.01e03 & 2.02e+03 \\ [2mm]
    & Intermediate & 6.08e-01 & 1.96e-01 & 1.51e-04 & 1.95e+03 & 9.41e-02 & 2.61e+01 & 6.12e+02 & 1.22e+03 \\ [2mm]
    & Fine & 4.45e-01 & 1.09e-01 & 1.60e-04 & 1.00e+03 & 4.84e-02 & 1.86e+01 & 3.23e+02 & 6.46e+02 \\
    \bottomrule
    \end{tabular}
    \end{adjustbox}
    \vspace{1ex}
    
    {\scriptsize\raggedright The results for $\mathcal{M}^{\oplus}_{\gothu;\mathrm{EC}}$ are omitted since they are equal to $\mathcal{M}^{\oplus}_{\gothp;\mathrm{EC}}$. \par}
\end{table}

In this numerical experiment, we consider the benchmark case~3b from~\cite{flemisch2018benchmarks}. As shown in the left panel of Figure~\ref{fig:md-examples}, the domain consists of ten (partially intersecting) fractures embedded in a unit square matrix. The exact fracture coordinates can be found in Appendix~C of~\cite{flemisch2018benchmarks}. Fractures 4 and 5 represent blocking fractures ($\mathcal{K} = 10^{-4}$ and $\kappa = 1$) whereas the others represent conductive fractures ($\mathcal{K} = 10^{4}$ and $\kappa = 10^{8}$). The matrix permeability is set to one. A linear pressure drop is imposed from left ($p=4$) to right ($p=1$), whereas no flux is prescribed at the top and bottom of the domain.

\begin{figure}
\centering
\includegraphics[width=\textwidth]{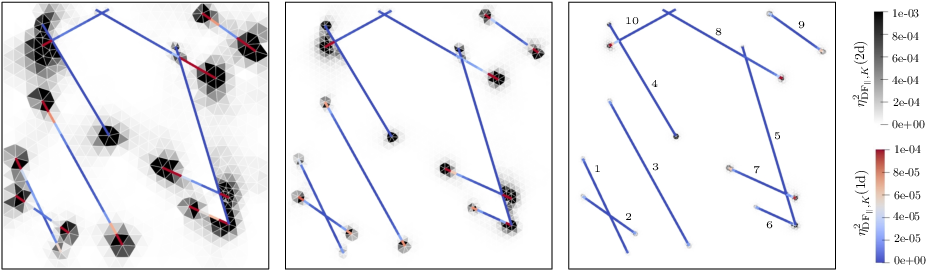}
\caption{Two-dimensional benchmark problem and the errors associated with the matrix and fractures for the coarse (left), intermediate (center), and fine (right) grid resolutions. Fractures 4 and 5 are blocking, whereas the others are conductive. The local bounds were obtained using MPFA. The results suggest that subdomain diffusive errors are concentrated around fracture tips and fracture intersections. \label{fig:bench2d}}
\end{figure}

The benchmark establishes three refinement levels; coarse, intermediate, and fine, with approximately $1500$, $4200$, and $16000$ two-dimensional cells. The structure of the local contributions to the majorant (confer e.g. equation \eqref{eq:df_loc}) are shown in Figure \ref{fig:bench2d}, based on the approximate solution obtained by the MPFA discretization.

In Table~\ref{tab:bench2d}, we show the errors bounds for the three refinement levels. To avoid numbering domains and interfaces, we refer to the matrix error as $\varepsilon_{\Omega^2,\mathrm{EC}}$, and group the fracture and interface errors by conductive and blocking. For example, $\varepsilon_{\Omega^1,\mathrm{C},\mathrm{EC}}$ refers to the sum of the errors of 1d conductive fractures.

An important observation is that the persistent reduction of the majorant $\mathcal{M}^{\oplus}_{\gothp,\gothu;\mathrm{EC}}$, together with the known upper and lower bounds on the efficiency indexes established in Remark \ref{rem:eff_ind}, provides a \textit{post factum} verification of the convergence of all the numerical methods considered. 

The error estimates suggest that the contribution to the overall error bounds are concentrated, primarily, on highly conductive interfaces (see the column corresponding to $\varepsilon_{\Gamma^1,\mathrm{C}}$). On a more qualitative note, Figure~\ref{fig:bench2d} suggests that subdomain diffusive errors are concentrated at the fracture tips and fracture intersections, which is where singularities may typically be encountered~\cite{boon2018robust}.

\subsection{Three-dimensional application\label{sec:3dbench}}

Our last numerical application is based on a modified version of the three-dimensional benchmark 2.1 from~\cite{flemisch2020verification}. The domain consists of nine intersecting fractures embedded in a unit cube, as shown in the middle panel of Figure~\ref{fig:md-examples}. This results in an intricate network with 106 subdomains and 270 interfaces of different dimensionality. 

\begin{table}
\caption{Error estimates for the three-dimensional application.\label{tab:geiger3Dscaled}}

    \begin{adjustbox}{width=\textwidth}
    \begin{tabular}{r r c c c c c c c c}
    \toprule
    & Mesh & ${{\varepsilon}}_{\Omega^3;\mathrm{EC}}$ & $ {{\varepsilon}}_{\Omega^2;\mathrm{EC}}$ & ${{\varepsilon}}_{\Omega^1;\mathrm{EC}}$ & $ {{\varepsilon}}_{\Gamma^2}$ & $ {{\varepsilon}}_{\Gamma^1}$ & $ {{\varepsilon}}_{\Gamma^0}$ & ${\mathcal{M}}^{\oplus}_{\gothp;\mathrm{EC}}$ &  ${\mathcal{M}}^{\oplus}_{\gothp,\gothu;\mathrm{EC}}$ \\
    \midrule 
    \multirow{3}{*}{\rotatebox[origin=c]{90}{RT0-P0~~}} & Coarse & 6.17e-01 & 5.81e-04 & 3.16e-04 & 9.87e+02 & 3.63e-02 & 3.31e-02 & 5.03e+02 & 1.01e+03 \\ [2mm]
    & Intermediate & 4.55e-01 & 4.61e-04 & 1.58e-04 & 7.75e+01 & 8.86e-03 & 8.35e-04 & 3.40e+01 & 6.81e+01 \\ [2mm]
    & Fine & 3.86e-01 & 2.55e-04 & 9.60e-05 & 2.26e+01 & 4.63e-03 & 4.34e-04 & 1.07e+01 & 2.14e+01 \\
    \midrule
    \multirow{3}{*}{\rotatebox[origin=c]{90}{MVEM-P0}} & Coarse & 6.07e-01 & 6.99e-04 & 2.77e-04 & 9.54e+02 & 7.48e-02 & 6.38e-02 & 4.66e+02 & 9.33e+02 \\ [2mm]
    & Intermediate & 4.55e-01 & 4.63e-04 & 1.65e-04 & 8.19e+01 & 9.96e-03 & 4.59e-03 & 3.59e+01 & 7.18e+01 \\ [2mm]
    & Fine & 3.86e-01 & 2.46e-04 & 9.17e-05 & 2.33e+01 & 4.00e-03 & 1.75e-03 & 1.11e+01 & 2.22e+01 \\
    \midrule
    \multirow{3}{*}{\rotatebox[origin=c]{90}{MPFA~~}} & Coarse & 6.07e-01 & 7.00e-04 & 3.15e-04 & 1.05e+03 & 4.61e-02 & 1.69e-02 & 5.24e+02 & 1.05e+03 \\ [2mm]
    & Intermediate & 4.46e-01 & 4.88e-04 & 1.61e-04 & 8.42e+01 & 7.72e-03 & 2.31e-03 & 3.71e+01 & 7.42e+01 \\ [2mm]
    & Fine & 3.77e-01 & 2.53e-04 & 9.04e-05 & 2.37e+01 & 2.82e-03 & 9.36e-04 & 1.12e+01 & 2.24e+01 \\
    \midrule
    \multirow{3}{*}{\rotatebox[origin=c]{90}{TPFA~~}} & Coarse & 6.32e-01 & 4.72e-04 & 2.26e-04 & 7.92e+02 & 4.21e-02 & 1.34e-02 & 3.76e+02 & 7.52e+02 \\ [2mm]
    & Intermediate & 4.48e-01 & 6.27e-04 & 1.40e-04 & 1.47e+02 & 1.56e-02 & 2.32e-03 & 6.82e+01 & 1.36e+02 \\ [2mm]
    & Fine & 4.07e-01 & 5.82e-04 & 8.72e-05 & 4.60e+01 & 7.97e-03 & 1.05e-03 & 2.04e+01 & 4.08e+01 \\
    \bottomrule
    \end{tabular}
    \end{adjustbox}
    \vspace{1ex}
    
    {\scriptsize\raggedright The results for $\mathcal{M}^{\oplus}_{\gothu;\mathrm{EC}}$ are omitted since they are equal to $\mathcal{M}^{\oplus}_{\gothp;\mathrm{EC}}$. \par}
\end{table}

\begin{figure}[tbp]
\centering
\includegraphics[height=0.25\textheight]{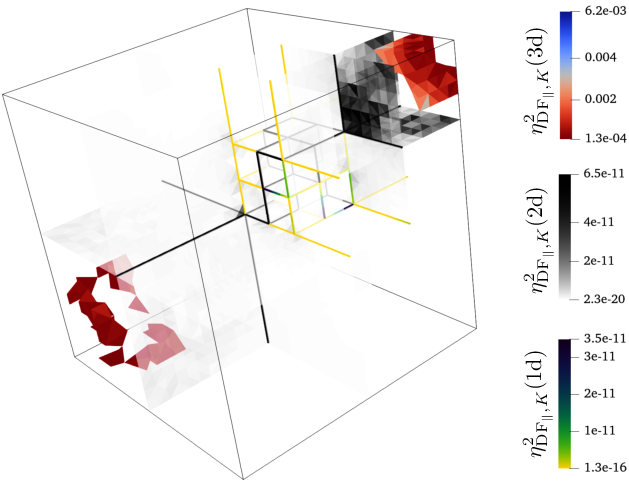}
\caption{Subdomain diffusive error contributions to the majorant for the fine grid resolution obtained with RT0-P0. \label{fig:geiger}}
\end{figure}

The original benchmark imposes an inlet flux (purple lower corner ${u}=-1$) and an outlet pressure (pink upper corner $p=1$), and for the rest of the external boundaries null flux. Since we have only detailed our results for zero Neumann boundary conditions, we have replaced the inlet flux by a constant pressure condition ($p=1$) and modified the value of the outlet pressure ($p=0$). The benchmark assigns heterogeneous permeability to the matrix subdomain, whereas the fractures are assumed to be highly conductive. For the complete description of the benchmark, we refer to~\cite{flemisch2020verification}, and for an impression on how the contributions to the majorant are distributed, see Figure~\ref{fig:geiger}. Here we show the error estimates for the whole fracture network obtained with RT0-P0, where it becomes evident that the subdomain diffusive errors are concentrated at the inlet and outlet boundaries; refinement efforts should therefore focus on these regions.

As in Section~\ref{sec:2dbench}, we collect the local errors of subdomains and interfaces of equal dimensionality. The results are summarized in Table~\ref{tab:geiger3Dscaled}. As in the previous cases, we have local and global convergence for all four numerical methods. Again, RT0-P0, MVEM-P0, and MPFA show very similar results, while TPFA show larger errors. 

As in the 2d case discussed above, the persistent reduction of the majorant $\mathcal{M}^{\oplus}_{\gothp,\gothu;\mathrm{EC}}$, again serves as a verification of the convergence of all four numerical methods. 

\section{Conclusion\label{sec:conclusion}}

In this paper, we obtained \textit{a posteriori} error estimates for mixed-dimensional elliptic equations. Depending upon the level of accuracy at which residual balances can be approximated, we have derived four concrete versions of the majorant; i.e.: for no mass-conservative, subdomain mass-conservative, locally mass-conservative, and point-wise mass-conservative approximations. Furthermore, we have demonstrated both theoretically and numerically that sharper bounds can be obtained (for locally mass-conservative methods) using local Poincar{\'e} constants instead of the global ones.

Our bounds have been thoroughly tested with numerical approximations obtained with four locally mass-conservative methods of the lowest-order, namely: RT0-P0, MVEM-P0, MPFA, and TPFA. We performed a detailed efficiency analysis comparing the use of global and local Poincar{\'e}-Friedrichs constants in two and three dimensions. In both validations, our upper bounds reflected the optimal convergence rates of the numerical methods. In addition, we applied our bounds to two- and three-dimensional community benchmark problems exhibiting challenging fracture networks. Again, in both cases, the bounds reflected the limitations and the convergence rates of the methods satisfactory.

To the best of our knowledge, the bounds obtained here are the first of their kind to provide a practical tool to measure the error in numerical approximations to the equations modeling the incompressible, single-phase flow in generic fractured porous media.

\vspace{2ex}\noindent\textbf{Funding:} \textit{Jhabriel Varela was funded by VISTA – a basic research program in collaboration between The Norwegian Academy of Science and Letters, and Equinor. Additionally, this work was supported in part through the Norwegian Research Council grant 250223. The authors would like to thank W. M. Boon for helpful discussions on this topic.}

\bibliographystyle{unsrtnat}
\bibliography{references}

\newpage
\begin{appendices}

\section{Derivation of variational formulations}

Here, we present the derivations for the primal and dual variational formulations for the case of a single fracture immersed in a matrix.

\subsection{Derivation of the primal weak form for a single fracture\label{sec:primalWeakFormProof}}

Substitute~\eqref{eq:DarcyLawMatrix} into~\eqref{eq:massConservationMatrix}, multiply each term by $q_2 \in H_{0}^1(\Omega_2)$, and integrate over $\Omega_2$. Similarly, substitute~\eqref{eq:DarcyLawMatrix}, \eqref{eq:DarcyLawInterface1}, and \eqref{eq:DarcyLawInterface2} into~\eqref{eq:massConservationFracture}, multiply each term by $q_1\in H_{0}^1(\Omega_1)$ and integrate over $\Omega_1$. Add the resulting equations to obtain
\begin{align}
   &-\inner{\nabla_2\cdot\mathcal{K}_2\nabla_2\,p_2}{q_2}_{\Omega_2} -\inner{\nabla_1\cdot\mathcal{K}_1\nabla_1\,p_1}{q_1}_{\Omega_1} + \inner{\kappa_{1}\left(p_1 - \mathrm{tr}_{\partial_1\Omega_2}\,p_2\right)}{q_1}_{\Omega_1} \nonumber\\
   &\quad + \inner{\kappa_{2}\left(p_1 - \mathrm{tr}_{\partial_2\Omega_2}\,p_2\right)}{q_1}_{\Omega_1}
   = \inner{f_2}{q_2}_{\Omega_2} + \inner{f_1}{q_1}_{\Omega_1}.\label{eq:primalWeakI}
\end{align}
   
Using integration by parts, the first term of~\eqref{eq:primalWeakI} can be expressed as
\begin{align}
    &-\inner{\nabla_2\cdot\mathcal{K}_2\nabla_2\,p_2}{q_2}_{\Omega_2} \nonumber \\
    &~~= \inner{\mathcal{K}_2\nabla_2\,p_2}{\nabla_2\,q_2}_{\Omega_2} - \sum_{j=1}^2 \innerangle{\mathrm{tr}_{\partial_j\Omega_2}\,(\mathcal{K}_2\nabla_2\,p_2)\cdot \vecn_{2}}{\mathrm{tr}_{\partial_j\Omega_2}\,q_2}_{\partial_j\Omega_2}, \nonumber\\
    &~~= \inner{\mathcal{K}_2\nabla_2p_2}{\nabla_2q_2}_{\Omega_2} - \sum_{j=1}^2 \inner{\lambda_j}{\mathrm{tr}_{\partial_j\Omega_2}\,q_2}_{\Gamma_j}, \nonumber\\
    &~~= \inner{\mathcal{K}_2\nabla_2p_2}{\nabla_2 q_2}_{\Omega_2} + \sum_{j=1}^2 \innerangle{\kappa_{j}\left(p_1 - \mathrm{tr}_{\partial_j\Omega_2}\,p_2\right)}{\mathrm{tr}_{\partial_j\Omega_2}\,q_2}_{\Gamma_{j}}.\label{eq:primalWeakII}
\end{align}
Here, we use the internal boundary conditions~\eqref{eq:InternalBoundary1} and \eqref{eq:InternalBoundary2} and the definition of the mortar fluxes~\eqref{eq:DarcyLawInterface1} and \eqref{eq:DarcyLawInterface2}. Analogously, integration by parts allows us to write the second term of~\eqref{eq:primalWeakI} as
\begin{equation}
    -\inner{\nabla_1\cdot \mathcal{K}_1 \nabla_1\, p_1}{q_1}_{\Omega_1} = \inner{\mathcal{K}_1 \nabla_1\,p_1}{\nabla_1\,q_1}_{\Omega_1}.\label{eq:primalWeakIII}
\end{equation}
Note that the boundary terms vanish due to the choice of boundary conditions. 

Finally, we note that the third and fourth terms from~\eqref{eq:primalWeakI} can be equivalently written as
\begin{equation}
    \inner{\kappa_{j}\left(p_1 - \mathrm{tr}_{\partial_j\Omega_2}\,p_2\right)}{q_1}_{\Omega_1} = \innerangle{\kappa_{j}\left(p_1 - \mathrm{tr}_{\partial_j\Omega_2}\,p_2\right)}{q_1}_{\Gamma_{j}}, \qquad j\in\{1,2\}.\label{eq:primalWeakIV}
\end{equation}
The proof is completed by substituting~\eqref{eq:primalWeakII}, \eqref{eq:primalWeakIII}, and \eqref{eq:primalWeakIV} into~\eqref{eq:primalWeakI} and grouping common terms.

\subsection{Derivation of the dual weak form for a single fracture\label{sec:dualWeakFormProof}}

Let us start with~\eqref{eq:dualFormVelocity}. Multiply respectively~\eqref{eq:DarcyLawMatrix} and~\eqref{eq:DarcyLawFracture} by $\vecv_{0,2}\in {V}_{0,2}$ and $\vecv_{0,1}\in {V}_{0,1}$, integrate over the subdomains $\Omega_2$ and $\Omega_1$, use integration by parts to obtain
\begin{align}
    \inner{\mathcal{K}_2^{-1}\vecu_{2}}{\vecv_{0,2}}_{\Omega_2} &= \inner{\mathcal{K}_2^{-1}\left(\vecu_{0,2} + \mathcal{R}_{1}\lambda_{1}+\mathcal{R}_{2}\lambda_{2}\right)}{\vecv_{0,2}}_{\Omega_2} = -\inner{\nabla_2\,p_2}{\vecv_{0,2}}_{\Omega_2} \nonumber \\
    & = \inner{p_2}{\nabla_2\cdot \vecv_{0,2}}_{\Omega_2} - \innerangle{g_{D,2}}{\mathrm{tr}_{\partial_D\Omega_2}\,\vecv_{0,2}\cdot\vecn_2}_{\partial_D\Omega_2}.\label{eq:dualProofVelocityI}\\
    \inner{\mathcal{K}_1^{-1}\vecu_{1}}{\vecv_{0,1}}_{\Omega_1} &= \inner{\mathcal{K}_1^{-1}\vecu_{0,1}}{\vecv_{0,1}}_{\Omega_1} = -\inner{\nabla_1\,p_1}{\vecv_{0,1}}_{\Omega_1} \nonumber \\
    & = \inner{p_1}{\nabla_1\cdot \vecv_{0,1}}_{\Omega_1} - \innerangle{g_{D,1}}{\mathrm{tr}_{\partial_D\Omega_1}\,\vecv_{0,1}\cdot\vecn_1}_{\partial_D\Omega_1}. \label{eq:dualProofVelocityII}
\end{align}
Adding together~\eqref{eq:dualProofVelocityI} and~\eqref{eq:dualProofVelocityII} gives~\eqref{eq:dualFormVelocity}. We now focus on~\eqref{eq:dualFormMortarFlux}. First, we use~\eqref{eq:DarcyLawMatrix} and multiply by the test functions $\mathcal{R}_j\nu_j$ with $\nu_{j}\in L^2(\Gamma_{j})$ for $j\in\{1,2\}$, integrate over $\Omega_2$, and apply integration by parts, to obtain:
\begin{align}
    \inner{\mathcal{K}_2^{-1}\vecu_{2}}{\mathcal{R}_{j}\nu_{j}}_{\Omega_2} &= \inner{\mathcal{K}_{2}^{-1}\left(\vecu_{0,2}+\mathcal{R}_{1}\lambda_{1}+\mathcal{R}_{2}\lambda_{2}\right)}{\mathcal{R}_{j}\nu_{j}}_{\Omega_2} =  -\inner{\nabla_2\,p_2}{\mathcal{R}_{j}\nu_{j}}_{\Omega_2} \nonumber \\
    &= \inner{p_2}{\nabla_2\cdot(\mathcal{R}_{j} \nu_{j})}_{\Omega_2} - \innerangle{\mathrm{tr}_{\partial_j\Omega_2}\,p_2}{\mathrm{tr}_{\partial_j\Omega_2}\,(\,\mathcal{R}_{j}\nu_{j})\cdot\vecn_2}_{\partial_j\Omega_2} \nonumber \\
    &= \inner{p_2}{\nabla_2\cdot(\mathcal{R}_{j} \nu_{j})}_{\Omega_2} - \innerangle{\mathrm{tr}_{\partial_j\Omega_2}\,p_2}{\nu_{j}}_{\partial_j\Omega_2} \nonumber \\
    &= \inner{p_2}{\nabla_2\cdot(\mathcal{R}_{j} \nu_{j})}_{\Omega_2} - \innerangle{\mathrm{tr}_{\partial_j\Omega_2}\,p_2}{\nu_{j}}_{\Gamma_{j}}. \label{eq:dualProofMortarFluxI}
\end{align}
Next, we multiply the interface laws~\eqref{eq:DarcyLawInterface1} and~\eqref{eq:DarcyLawInterface2} by $\nu_{1}$ and $\nu_2$, respectively, to get for $j=\{1,2\}$
\begin{equation}
    \innerangle{\kappa_{1}^{-1}\lambda_{j}}{\nu_{j}}_{\Gamma_{j}} = - \innerangle{p_1}{\nu_{j}}_{\Gamma_{j}} + \innerangle{\mathrm{tr}_{\partial_j\Omega_2}\,p_2}{\nu_{j}}_{\Gamma_{j}} = -\inner{p_1}{\nu_{j}}_{\Omega_1} + \innerangle{\mathrm{tr}_{\partial_j\Omega_2}\,p_2}{\nu_{j}}_{\Gamma_{j}}. \label{eq:dualProofMortarFluxII}
\end{equation}
After adding~\eqref{eq:dualProofMortarFluxI} and \eqref{eq:dualProofMortarFluxII} and canceling common terms, we obtain~\eqref{eq:dualFormMortarFlux}. Finally, to obtain~\eqref{eq:dualFormPressure}, we multiply~\eqref{eq:massConservationMatrix} by $q_2\in L^2(\Omega_2)$ and ~\eqref{eq:massConservationFracture} by $q_1\in L^2(\Omega_1)$, and integrate over their respective subdomains, and add the resulting equations.

\section{Proof of Theorem~\ref{thm:firstUpperBound}\label{sec:firstThmProof}}

Here, we present the proof of the upper bound of the error for the primal variable, for the case of a single fracture immersed in a matrix.

\begin{proof}
Start by computing the difference between $p=[p_1,p_2]\in H_{0}^1(\Omega)+g$ and an arbitrary function $q=[q_1,q_2]\in H_{0}^1(\Omega)+g$ in the energy norm~\eqref{eq:firstEnergyNorm}:
\begin{align}
\tnorm{p-q}^2 &= \inner{\mathcal{K}_2\nabla_2(p_2-q_2)}{\nabla_2(p_2-q_2)}_{\Omega_2} + \inner{\mathcal{K}_1\nabla_1(p_1-q_1)}{\nabla_1(p_1-q_1)}_{\Omega_1} \nonumber \\
&~~~+\sum_{j=1}^2\innerangle{\kappa_{j}\left[(p_1-q_1)-\mathrm{tr}_{\partial_j\Omega_2}\,(p_2-q_2)\right]}{(p_1-q_1)-\mathrm{tr}_{\partial_j\Omega_2}\,(p_2-q_2)}_{\Gamma_{j}}, \nonumber \\
&= \inner{\mathcal{K}_2\nabla_2p_2}{\nabla_2(p_2-q_2)}_{\Omega_2} + \inner{\mathcal{K}_1\nabla_1p_1}{\nabla_1(p_1-q_1)}_{\Omega_1} \nonumber \\
&~~~+ \sum_{j=1}^2\innerangle{\kappa_{j}\left[(p_1-q_1)-\mathrm{tr}_{\partial_j\Omega_2}(p_2-q_2)\right]}{(p_1-q_1)-\mathrm{tr}_{\partial_j\Omega_2}\,(p_2-q_2)}_{\Gamma_{j}} \nonumber \\
&~~~+\inner{-\mathcal{K}_2\nabla_2q_2}{\nabla_2(p_2-q_2)}_{\Omega_2} + \inner{-\mathcal{K}_1\nabla_1q_1}{\nabla_1(p_1-q_1)}_{\Omega_1} \nonumber \\
&~~~+ \sum_{j=1}^2\innerangle{-\kappa_{j}\left(q_1-\mathrm{tr}_{\partial_j\Omega_2}\,q_2\right)}{(p_1-q_1)-\mathrm{tr}_{\partial_j\Omega_2}\,(p_2-q_2)}_{\Gamma_{j}}.\label{eq:fproof1}
\end{align}
By noticing that the first three terms of~\eqref{eq:fproof1} add up to the right-hand side of~\eqref{eq:primalSingleFracture}, and adding the identity
\begin{align}
    &-\inner{\vecv_{0,2}+\mathcal{R}_{1}\nu_{1}+\mathcal{R}_{2}\nu_{2}}{\nabla_2(p_2-q_2)}_{\Omega_2} - \inner{\vecv_{0,1}}{\nabla_1(p_1-q_1)}_{\Omega_1} \nonumber\\
    &~~~+ \inner{\nabla_2\cdot\left(\vecv_{0,2}+\mathcal{R}_{1}\nu_{1}+\mathcal{R}_{2}\nu_{2}\right)}{p_2-q_2}_{\Omega_2} + \inner{\nabla_1\cdot\vecv_{0,1}-\nu_1-\nu_2}{p_1-q_1}_{\Omega_1} \nonumber \\
    &~~~+ \sum_{j=1}^2 \innerangle{\nu_{j}}{(p_1-q_1) - \mathrm{tr}_{\partial_j\Omega_2}\,(p_2-q_2)}_{\Gamma_{j}} = 0,\nonumber
\end{align}
valid for any $\vecv_{0}\in{V}_0$ and $\nu\in L^2(\Gamma)$ to \eqref{eq:fproof1}, we obtain
\begin{align}
    \tnorm{p-q}^2 &= \inner{-\left(\vecv_{0,2}+\mathcal{R}_{1}\nu_{1}+\mathcal{R}_{2}\nu_{2}+\mathcal{K}_2\nabla_2q_2\right)}{\nabla_2(p_2-q_2)}_{\Omega_2} \nonumber \\
    &~~~~~~+ \inner{-\left(\vecv_{0,1}+\mathcal{K}_1\nabla_1p_1\right)}{\nabla_1(p_1-q_1)}_{\Omega_1} \nonumber \\
    &~~~~~~+ \sum_{j=1}^2\innerangle{-\left[\nu_{j}+\kappa_{j}\left(q_1-\mathrm{tr}_{\partial_j\Omega_2}\,q_2\right)\right]}{(p_1-q_1)-\mathrm{tr}_{\partial_j\Omega_2}\,(p_2-q_2)}_{\Gamma_{j}} \nonumber \\
    &~~~~~~+ \inner{f_2-\nabla_2\cdot\left(\vecv_{0,2}+\mathcal{R}_{1}\nu_{1}+\mathcal{R}_{2}\nu_{2}\right)}{p_2-q_2}_{\Omega_2} \nonumber \\
    &~~~~~~+ \inner{f_1-\nabla_1\cdot\vecv_{0,1}+\nu_{1}+\nu_2}{p_1-q_1}_{\Omega_1}.\label{eq:fproof2}
\end{align}

Recognizing that since $\mathcal{K}_2$ is symmetric positive definite, it can be expressed as $\mathcal{K}_2 = \left(\mathcal{K}_2^{1/2}\right)^2$, where $\mathcal{K}_2^{1/2}$ is also symmetric positive definite, and therefore self-adjoint. The square-root of the material coefficients can therefore be moved to the second argument of the three first inner products in \eqref{eq:fproof2}. After applying the Cauchy-Schwarz inequality to each inner product of \eqref{eq:fproof2}, one gets
\begin{align}
    &\tnorm{p-q}^2 \leq \norm{\mathcal{K}_2^{\mhalf}\left(\vecv_{0,2}+\mathcal{R}_{1}\nu_{1}+\mathcal{R}_{2}\nu_{2}+\mathcal{K}_2\nabla_2q_2\right)}_{\Omega_2} ~\norm{\mathcal{K}_2^\phalf\nabla_2(p_2-q_2)}_{\Omega_2} \nonumber \\
    &+ \norm{\mathcal{K}_1^{\mhalf}\left(\vecv_{0,1}+\mathcal{K}_1\nabla_1p_1\right)}_{\Omega_1} ~\norm{\mathcal{K}_1^{\phalf}\nabla_1(p_1-q_1)}_{\Omega_1} \nonumber\\
    &+ \norm{\kappa_{1}^{\mhalf}\left[\nu_{1}+\kappa_{1}\left(q_1-\mathrm{tr}_{\partial_1\Omega_2}\,q_2\right)\right]}_{\Gamma_{1}}~ \norm{\kappa_{1}^{\phalf}\left[(p_1-q_1)-\mathrm{tr}_{\partial_1\Omega_2}\,(p_2-q_2)\right]}_{\Gamma_{1}} \nonumber \\
    &+ \norm{\kappa_{2}^{\mhalf}\left[\nu_{2}+\kappa_{2}\left(q_1-\mathrm{tr}_{\partial_2\Omega_2}\,q_2\right)\right]}_{\Gamma_{2}}~ \norm{\kappa_{2}^{\phalf}\left[(p_1-q_1)-\mathrm{tr}_{\partial_2\Omega_2}\,(p_2-q_2)\right]}_{\Gamma_{2}} \nonumber \\
    &+ \norm{f_2-\nabla_2\cdot\left(\vecv_{0,2}+\mathcal{R}_{1}\nu_{1}+\mathcal{R}_{2}\nu_{2}\right)}_{\Omega_2} ~\norm{p_2-q_2}_{\Omega_2} \nonumber \\
    &+ \norm{f_1-\nabla_1\cdot\vecv_{0,1}+\nu_1 +\nu_2}_{\Omega_1} ~\norm{p_1-q_1}_{\Omega_1}\nonumber
\end{align}
Applying the permeability-weighted Poincar{\'e}-Friedrichs inequality~\eqref{eq:CP_subdomain} to the terms $\norm{p_1-q_1}_{\Omega_1}$ and $\norm{p_2-q_2}_{\Omega_2}$, the proof of the theorem is completed.
\end{proof}

\section{Proof of Theorem~\ref{thm:Abstract}\label{sec:thmProof}}

Here, we present the proof of our main theorem, which deals with the general abstract estimates in a mixed-dimensional setting.

\begin{proof}
(1) The proof for the bounds for the mD primal variable follows the one presented in Appendix~\ref{sec:firstThmProof}, modulo its generalization to the mD~setting and the use of weighted norms on the residual terms. Start by computing the difference between any $\gothq \in \mixHIzero+\gothg$ and $\gothp\in\mixHIzero+\gothg$ using~\eqref{eq:bilinearB}, to get
\begin{align}
\tnorm{\gothp - \gothq}^2 &= \inner{\mixperm\mixgrad(\gothp-\gothq)}{\mixgrad(\gothp-\gothq)}_{\Omega, \Gamma} \nonumber \\
&= \inner{\mixperm\mixgrad\gothp}{\mixgrad(\gothp-\gothq)}_{\Omega, \Gamma} + \inner{-\mixperm\mixgrad\gothq}{\mixgrad(\gothp-\gothq)}_{\Omega, \Gamma}\nonumber \\
&= \inner{\gothf}{\gothp-\gothq}_{\Omega} + \inner{-\mixperm\mixgrad\gothq}{\mixgrad(\gothp-\gothq)}_{\Omega, \Gamma} \nonumber \\
&= \inner{\gothf}{\gothp-\gothq}_{\Omega} + \inner{-\mixperm^{\mhalf}\mixgrad\gothq}{\mixperm^{\phalf}\mixgrad(\gothp-\gothq)}_{\Omega, \Gamma} \nonumber \\
&= \inner{\gothf-\mixdiv\gothv}{\gothp-\gothq}_{\Omega} + \inner{-\mixperm^{\mhalf}(\gothv + \mixperm\mixgrad\gothq)}{\mixperm^{\phalf}\mixgrad(\gothp-\gothq)}_{\Omega, \Gamma}. \label{eq:proofPrimal_a}
\end{align}
Here, we used \eqref{eq:bilinearB}, \eqref{eq:primalAbstract}, and added the fact that $\mixdiv$ and $\mixgrad$ are adjoints. 

By exploiting the orthogonality property~\eqref{eq:ortho} and then introducing the weights to the second and third terms, \eqref{eq:proofPrimal_a} can be equivalently written as:
\begin{align}
 &\tnorm{\gothp - \gothq}^2 = \inner{\gothf-\mixdiv\gothv}{\pi_{W}(\gothp-\gothq)}_{\Omega} +  \inner{-\mixperm^{\mhalf}(\gothv + \mixperm\mixgrad\gothq)}{\mixperm^{\phalf}\mixgrad(\gothp-\gothq))}_{\Omega, \Gamma} \nonumber \\
&= \inner{\mu^{-1}(\gothf-\mixdiv\gothv)}{\mu\pi_{W}(\gothp-\gothq)}_{\Omega} + \inner{-\mixperm^{\mhalf}(\gothv + \mixperm\mixgrad\gothq)}{\mixperm^{\phalf}\mixgrad(\gothp-\gothq)}_{\Omega, \Gamma}. \label{eq:proofPrimal_b}
\end{align}

Finally, applying the Cauchy-Schwarz inequality to the first and second terms of~\eqref{eq:proofPrimal_b}, and then the norm definitions~\eqref{eq:bilinearB}, \eqref{eq:bilinearA}, and \eqref{eq:weightedNorm}, we arrive at the desired bound:
\begin{align}
   &\tnorm{\gothp - \gothq}^2  \leq \tnormstar{\gothv+\mixperm\mixgrad\gothq} \tnorm{\gothp-\gothq} + \norm{\mu^{-1}(\gothf-\mixdiv\gothv)}_{\Omega} \norm{\pi_{W}(\gothp-\gothq)}_{W,\mu} \nonumber \\ 
&\leq \tnormstar{\gothv+\mixperm\mixgrad\gothq} \tnorm{\gothp-\gothq} + \norm{\mu^{-1}(\gothf-\mixdiv\gothv)}_{\Omega} \tnorm{\gothp-\gothq}\leq \mathcal{M}(\gothq,\gothv,\gothf,\mu) \tnorm{\gothp-\gothq}.
\end{align}

(2) The proof for the bounds for the dual variable is given next. We remark that an alternative proof based on a generalized abstract estimate (see~\cite{vohralik2010unified}, Theorem~6.1) can be used to obtain equivalent upper bounds after its generalization to the mD~setting. 

We start by adding the square of the primal and dual error to obtain:
\begin{align}
    &\tnorm{\gothp-\gothq}^2 + \tnormstar{\gothu-\gothv}^2 = \inner{\mixperm\mixgrad(\gothp-\gothq)}{\mixgrad(\gothp-\gothq)}_{\Omega, \Gamma} + \inner{\mixperm^{-1}(\gothu-\gothv)}{\gothu-\gothv}_{\Omega, \Gamma} \nonumber\\
    &= \inner{\gothu + \mixperm\mixgrad\gothq}{\mixperm^{-1}\gothu+\mixgrad\gothq}_{\Omega, \Gamma} + \inner{\mixperm^{-1}(\gothu-\gothv)}{\gothu-\gothv}_{\Omega, \Gamma} \nonumber\\
    &= \inner{\gothu -\gothv + \gothv + \mixperm\mixgrad\gothq}{\mixperm^{-1}\gothu-\mixperm^{-1}\gothv+\mixgrad\gothq + \mixperm^{-1}\gothv}_{\Omega, \Gamma} + \inner{\mixperm^{-1}(\gothu-\gothv)}{\gothu-\gothv}_{\Omega, \Gamma} \nonumber\\
    &=\inner{\gothv+\mixperm\mixgrad\gothq}{\mixperm^{-1}\gothv+\mixgrad\gothq}_{\Omega, \Gamma} + 2\inner{\gothu-\gothv}{-\mixgrad(\gothp-\gothq)}_{\Omega, \Gamma} \nonumber\\
    &=\inner{\mixperm^{\mhalf}\gothv+\mixperm^{\phalf}\mixgrad\gothq}{\mixperm^{\mhalf}\gothv+\mixperm^\phalf\mixgrad\gothq}_{\Omega, \Gamma} + 2\inner{\gothu-\gothv}{-\mixgrad(\gothp-\gothq)}_{\Omega, \Gamma}.\label{eq:dualProof_a}
\end{align}

Here, we used the norm definitions~\eqref{eq:bilinearB} and~\eqref{eq:bilinearA} together with the mD~constitutive relationship~\eqref{eq:mD-Darcylaw}.

Using partial integration, mass conservation~\eqref{eq:mD-massConservation}, and the orthogonality property~\eqref{eq:ortho}, the second term of~\eqref{eq:dualProof_a} can be equivalently written as
\begin{align}
    &\inner{\gothu-\gothv}{-\mixgrad(\gothp-\gothq)}_{\Omega, \Gamma} = \inner{\mixdiv(\gothu-\gothv)}{-(\gothp-\gothq)}_\Omega = \inner{\gothf-\mixdiv\gothv}{-(\gothp-\gothq)}_\Omega \nonumber\\
    &= \inner{\gothf-\mixdiv\gothv}{-\pi_{W}(\gothp-\gothq)}_\Omega = \inner{\mu^{-1}(\gothf-\mixdiv\gothv)}{-\mu\pi_{W}(\gothp-\gothq)}_\Omega.\label{eq:dualProof_b}
\end{align}

Using the Cauchy-Schwarz inequality twice and the definition of the weighted norms~\eqref{eq:weightedNorm},~\eqref{eq:dualProof_b}
can be estimated as
\begin{align}
    &\abs{\inner{\gothu-\gothv}{-\mixgrad(\gothp-\gothq)}_{\Omega, \Gamma}} \leq \norm{\mu^{-1}(\gothf-\mixdiv\gothv)}_\Omega \norm{\pi_{W}(\gothp-\gothq)}_{W,\mu} \nonumber\\
    &~~~=\norm{\mu^{-1}(\gothf-\mixdiv\gothv)}_\Omega \tnorm{\gothp-\gothq} \leq \frac{1}{2} \left(\norm{\mu^{-1}(\gothf-\mixdiv\gothv)}_\Omega^2 +  \tnorm{\gothp-\gothq}^2\right).\label{eq:dualProof_c}
\end{align}

Substituting~\eqref{eq:dualProof_c} into~\eqref{eq:dualProof_a} and applying the Cauchy-Schwarz inequality to the first term, we arrive at,
\begin{align*}
    \tnormstar{\gothu-\gothv}^2 &\leq \tnormstar{\gothv+\mixperm\mixgrad\gothq}^2 +  \norm{\mu^{-1}(\gothf-\mixdiv\gothv)}_\Omega^2,
\end{align*}
from which we conclude that~\eqref{eq:dualAbstract} indeed holds.

(3) To prove the upper bound for the primal-dual pair, we choose an arbitrary pair $(\gothq,\gothv)\in (\mixHIzero+\gothg)\times\mixHdivU$, and measure its difference with the exact solution $(\gothp,\gothu) \in (\mixHIzero+\gothg)\times\mixHdiv$ in the norm~\eqref{eq:combinedNorm}, to get
\begin{align*}
    \fullnorm{(\gothp-\gothq,\gothu-\gothv)} &= \tnorm{\gothp-\gothq} + \tnormstar{\gothu-\gothv} + \norm{\mu^{-1}\mixdiv (\gothu-\gothv)}_{\Omega} \\
    &\leq 2\mathcal{M} + \norm{\mu^{-1}\mixdiv (\gothu-\gothv)}_{\Omega},
\end{align*}

where we use the bounds~\eqref{eq:primalAbstract} and \eqref{eq:dualAbstract}. 

For the proof of the lower bound, we start from the definition of the majorant, to get
\begin{align*}
    \mathcal{M} &= \tnormstar{\gothv+\mixperm\mixgrad\gothq} + \norm{\mu^{-1}(\gothf - \mixdiv \gothv)}_{\Omega} \\
    &\leq \tnormstar{\gothu-\gothv} + \tnormstar{\mixperm\mixgrad(\gothp-\gothq)} + \norm{\mu^{-1}(\gothf - \mixdiv \gothv)}_{\Omega} = \fullnorm{(\gothp-\gothq,\gothu-\gothv)}.
\end{align*}

This completes the proof for the two-sided bounds and the abstract theorem.
\end{proof}

\section{Exact solutions to numerical validations\label{sec:exactSolutions}}

Herein, we provide the exact expressions for the pressure, velocities, mortar fluxes, and source terms for the numerical validations presented in Section \ref{sec:num_val}.

We will conveniently define the following quantities for notational compactness:
\begin{align*}
    \alpha(\vecx) &= x_1 - 0.50, \\
    \beta_1(\vecx) &= x_2 - 0.25, ~~~~~~ \beta_2(\vecx) = x_2 - 0.75, \\
    \gamma_1(\vecx) &= x_3 - 0.25, ~~~~~~ \gamma_2(\vecx) = x_3 - 0.75,
\end{align*}
where $x = [x_1, x_2, x_3]$.

\subsection{Exact solutions for the 1d/2d validation}

The matrix subdomain $\Omega_2$ is decomposed into three regions, i.e. $\Omega_2 = \cup_{k=1}^3 \Omega_2^k$, given by:
\begin{align*}
    \Omega_2^1 &= \left\lbrace\vecx\in\Omega_2 : 0.00 < x_2 < 0.25 \right\rbrace, \\
    \Omega_2^2 &= \left\lbrace\vecx\in\Omega_2 : 0.25 \leq x_2 < 0.75 \right\rbrace, \\
    \Omega_2^3 &= \left\lbrace\vecx\in\Omega_2 : 0.75 \leq x_2 < 1.00 \right\rbrace.
\end{align*}

Let us now define the distance function $d(\vecx)$ from $\Omega_2$ to $\Omega_1$. That is,
\begin{equation}
    d(\vecx) = 
    \begin{cases}
        \left(\alpha(\vecx)^2 + \beta_1(\vecx)^2\right)^{0.5}, & \vecx\in\Omega_2^1 \\
        \left(\alpha(\vecx)^2\right)^{0.5}, & \vecx\in\Omega_2^2 \\
        \left(\alpha(\vecx)^2 + \beta_2(\vecx)^2\right)^{0.5}, & \vecx\in\Omega_2^3 \\
    \end{cases},
\end{equation}
and the bubble function $\omega(\vecx)$:
\begin{equation}
    \omega(\vecx) = 
    \begin{cases}
        \beta_1(\vecx)^2 \beta_2(\vecx)^2, & \vecx\in\Omega_2^2 \\
        0, & \mathrm{otherwise}
    \end{cases}.
\end{equation}

In Table~\ref{tab:exactSolution2D}, we include the exact solutions for all the variables of interest. Note that the parameter $n$ controls the regularity of the solution. For this particular validation, a value of $n=1.5$ was adopted.
\begin{table}
    \caption{Exact solutions for the 1d/2d validation.\label{tab:exactSolution2D}}
    \begin{adjustbox}{width=\textwidth}
    \begin{tabular}{r l l}
    \toprule
         \multirow{2}{*}{$p_2=$} & $d^{n+1} + \omega d$ & $\Omega_2^2$\\[3pt]
          & $d_2^{n+1}$ & $\Omega_2\setminus\Omega_2^2$ \\
          \midrule
          \multirow{3}{*}{\vspace{-12.5pt}$\vecu_2=$} & $-d^{n+1}(n+1)\left[\begin{matrix} \alpha & \beta_1\end{matrix}\right]$ & $\Omega_2^1$\\[5pt]
          & $-d\left[\begin{matrix}\alpha^{-1}\left(\omega + d^n(n+1)\right) &  2\beta_1^2\beta_2+2\beta_1\beta_2^2\end{matrix}\right]$ & $\Omega_2^2$\\[3pt]
          & $-d^{n+1}(n+1)\left[\begin{matrix} \alpha & \beta_2\end{matrix}\right]$ & $\Omega_2^3$\\[3pt]
          \midrule
          \multirow{3}{*}{\vspace{-6pt}$f_2=$} & $-d^{-2}(n+1)\left( 2d^{n+1} + \alpha^2d^{n-1}(n-1) + \beta_1^2d^{n-1}(n-1)\right)$ & $\Omega_2^1$\\[3pt]
          & $-2d\left(\beta_1(\beta_1+2\beta_2) + \beta_2(2\beta_1+\beta_2)
         \right) - d^{n-1} n (n + 1)$ & $\Omega_2^2$\\[3pt]
          & $-d^{-2}(n+1)\left(2d^{n+1} +\alpha^2d^{n-1}(n-1) + \beta_2^2d^{n-1}(n-1)\right)$ & $\Omega_2^3$\\
          \midrule
          $\lambda_{1}=$ & $\omega$ & $\Gamma_{1}$ \\
          $\lambda_{2}=$ & $\omega$ & $\Gamma_{2}$ \\
          \midrule
          $p_2 =$ & $0$ & $\partial_1\Omega_2$ \\
          $p_2 =$ & $0$ & $\partial_2\Omega_2$ \\
          \midrule
          $p_1=$ & $-\omega$ & $\Omega_1$ \\
          \midrule
          $\vecu_1=$ & $\left[\begin{matrix}0 &  2\beta_1 ^ 2  \beta_2 + 2\beta_1  \beta_2 ^ 2\end{matrix}\right]$ & $\Omega_1$ \\
          \midrule
          $\sum_{j\in\hat{S}_1}\lambda_j=$ & $2\omega$ & $\Omega_1$ \\
          \midrule
          $f_1=$ & $8\beta_1\beta_2 + 2(\beta_1^2 + \beta_2^2) - 2\omega$  & $\Omega_1$ \\
          \bottomrule
    \end{tabular}
    \end{adjustbox}
\end{table}

\subsection{Exact solutions for the 2d/3d validation}

Analogously to the previous case, we decompose the three-dimensional matrix $\Omega_2$ into nine subdomains, i.e. $\Omega_2 = \cup_{k=1}^9 \Omega_2^k$, given by
\begin{align*}
    \Omega_2^1 &= \left\lbrace \vecx\in\Omega_2 : 0.00 < x_2 < 0.25,\,\, 0.00 < x_3 < 0.25\right\rbrace, \\
    \Omega_2^2 &= \left\lbrace \vecx\in\Omega_2 : 0.00 < x_2 < 0.25,\,\, 0.25 \leq x_3 < 0.75\right\rbrace, \\
    \Omega_2^3 &= \left\lbrace \vecx\in\Omega_2 : 0.00 < x_2 < 0.25,\,\, 0.75 \leq x_3 < 1.00\right\rbrace, \\
    \Omega_2^4 &= \left\lbrace \vecx\in\Omega_2 : 0.25 \leq x_2 < 0.75,\,\, 0.00 < x_3 < 0.25\right\rbrace, \\
    \Omega_2^5 &= \left\lbrace \vecx\in\Omega_2 : 0.25 \leq x_2 < 0.75,\,\, 0.25 \leq x_3 < 0.75\right\rbrace, \\
    \Omega_2^6 &= \left\lbrace \vecx\in\Omega_2 : 0.25 \leq x_2 < 0.75,\,\, 0.75\leq x_3 < 1.00\right\rbrace, \\
    \Omega_2^7 &= \left\lbrace \vecx\in\Omega_2 : 0.75 \leq x_2 < 1.00,\,\, 0.00 < x_3 < 0.25\right\rbrace, \\
    \Omega_2^8 &= \left\lbrace \vecx\in\Omega_2 : 0.75 \leq x_2 < 1.00,\,\, 0.25 \leq x_3 < 0.75\right\rbrace, \\
    \Omega_2^9 &= \left\lbrace \vecx\in\Omega_2 : 0.75 \leq x_2 < 1.00,\,\, 0.75 \leq x_3 < 1.00\right\rbrace.
\end{align*}

The distance function $d_2(\vecx)$ from $\Omega_2$ to $\Omega_1$ is now given by
\begin{equation}
    d_2(\vecx) = 
        \begin{cases}
        \left(\alpha(\vecx)^2 + \beta_1(\vecx)^2 + \gamma_1(\vecx)^2\right)^{0.5}, & \vecx\in\Omega_2^1, \\
        \left(\alpha(\vecx)^2 + \beta_1(\vecx)^2\right)^{0.5}, & \vecx\in\Omega_2^2, \\
        \left(\alpha(\vecx)^2 + \beta_1(\vecx)^2 + \gamma_2(\vecx)^2\right)^{0.5}, & \vecx\in\Omega_2^3, \\
        \left(\alpha(\vecx)^2 + \gamma_1(\vecx)^2\right)^{0.5}, & \vecx\in\Omega_2^4, \\
        \left(\alpha(\vecx)^2\right)^{0.5}, & \vecx\in\Omega_2^5, \\
        \left(\alpha(\vecx)^2 + \gamma_2(\vecx)^2\right)^{0.5}, & \vecx\in\Omega_2^6, \\
        \left(\alpha(\vecx)^2 + \beta_2(\vecx)^2 + \gamma_1(\vecx)^2\right)^{0.5}, & \vecx\in\Omega_2^7, \\
        \left(\alpha(\vecx)^2 + \beta_2(\vecx)^2\right)^{0.5}, & \vecx\in\Omega_2^8, \\
        \left(\alpha(\vecx)^2 + \beta_2(\vecx)^2 + \gamma_2(\vecx)^2\right)^{0.5}, & \vecx\in\Omega_2^9,
    \end{cases}
\end{equation}
and the bubble function $\omega(\vecx)$:
\begin{equation}
    \omega(\vecx) = 
    \begin{cases}
        \beta_1(\vecx)^2 \beta_2(\vecx)^2 \gamma_1(\vecx)^2 \gamma_2(\vecx)^2, & \vecx\in\Omega_2^5 \\
        0, & \mathrm{otherwise}
    \end{cases}.
\end{equation}

In Table~\ref{tab:exactSolution3D}, we show the exact solutions for all the variables of interest. Once again, a value of $n = 1.5$ is adopted for this validation.
\begin{table}[!h]
    \caption{Exact solutions for the $2d/3d$ validation.\label{tab:exactSolution3D}}
    \begin{adjustbox}{width=\textwidth}
    \begin{tabular}{r l l}
         \toprule
         \multirow{2}{*}{$p_2=$} & $d^{n+1} + \omega d$ & $\Omega_2^2$\\
          & $d_2^{n+1}$ & $\Omega_2\setminus\Omega_2^2$ \\
          \midrule
          \multirow{9}{*}{$\vspace{-35pt}\vecu_2=$} & $-d^{n-1}(n+1)\left[\begin{matrix}\alpha & \beta_1 & \gamma_1\end{matrix}\right]$ & $\Omega_2^1$\\[3pt]
          & $-d^{n-1}(n+1)\left[\begin{matrix}\alpha & \beta_1 & 0\end{matrix}\right]$ & $\Omega_2^2$ \\[3pt]
          & $-d^{n-1}(n+1)\left[\begin{matrix}\alpha & \beta_1 & \gamma_2\end{matrix}\right]$ & $\Omega_2^3$\\[3pt]
          &$-d^{n-1}(n+1)\left[\begin{matrix}\alpha & 0 & \gamma_1\end{matrix}\right]$ & $\Omega_2^4$\\[3pt]
          &$-d\left[\begin{matrix}\alpha^{-1}(\omega + d^n(n+1)) & 2 \beta_{1}^{2} \beta_{2} \gamma_{1}^{2} \gamma_{2}^{2} + 2 \beta_{1} \beta_{2}^{2} \gamma_{1}^{2} \gamma_{2}^{2} & 2 \beta_{1}^{2} \beta_{2}^{2} \gamma_{1}^{2} \gamma_{2} + 2 \beta_{1}^{2} \beta_{2}^{2} \gamma_{1} \gamma_{2}^{2}\end{matrix}\right]$ & $\Omega_2^5$ \\[3pt]
          &$-d^{n-1}(n+1)\left[\begin{matrix}\alpha & 0 & \gamma_2\end{matrix}\right]$ & $\Omega_2^6$\\[3pt]
          &$-d^{n-1}(n+1)\left[\begin{matrix}\alpha & \beta_2 & \gamma_1\end{matrix}\right]$ & $\Omega_2^7$\\[3pt]
          &$-d^{n-1}(n+1)\left[\begin{matrix}\alpha & \beta_2 & 0\end{matrix}\right]$ & $\Omega_2^8$\\[3pt]
          &$-d^{n-1}(n+1)\left[\begin{matrix}\alpha & \beta_2 & \gamma_2\end{matrix}\right]$ & $\Omega_2^9$\\[3pt]
          \midrule
          \multirow{10}{*}{$\vspace{-35pt}f_2=$} & $-d^{-2}(n+1) \left(3d^{n+1} + \alpha^2d^{n-1}(n-1) + \beta_1^2d^{n-1}(n-1) + \gamma_1^2d^{n-1}(n-1)\right)$ & $\Omega_2^1$\\[3pt]
          & $-d^{-2}(n+1) \left(2d^{n+1} + \alpha^2d^{n-1}(n-1) + \beta_1^2d^{n-1}(n-1) \right)$ & $\Omega_2^2$\\[3pt]
          &$-d^{-2}(n+1) \left(3d^{n+1} + \alpha^2d^{n-1}(n-1) + \beta_1^2d^{n-1}(n-1) + \gamma_2^2d^{n-1}(n-1)\right)$ & $\Omega_2^3$\\[3pt]
          &$-d^{-2}(n+1) \left(2d^{n+1} + \alpha^2d^{n-1}(n-1) + \gamma_1^2d^{n-1}(n-1) \right)$ & $\Omega_2^4$\\[3pt]
          &$-2d\left(\beta_{1}^{2} \beta_{2}^{2} \left(\gamma_{1} \left(\gamma_{1} + 2 \gamma_{2}\right) + \gamma_{2} \left(2 \gamma_{1} + \gamma_{2}\right)\right) +  \gamma_{1}^{2} \gamma_{2}^{2} \left(\beta_{1} \left(\beta_{1} + 2 \beta_{2}\right) + \beta_{2} \left(2 \beta_{1} + \beta_{2}\right)\right)\right)$& \multirow{2}{*}{$\Omega_2^5$}\\[5pt]
          & $\qquad\qquad -\alpha^{-2}\omega d^{n+1}(n+1)^2 - \alpha^{-2}\omega d^{n+1}(n+1)$ & \\[3pt]
          &$-d^{-2}(n+1) \left(2d^{n+1} + \alpha^2d^{n-1}(n-1) + \gamma_2^2d^{n-1}(n-1)\right)$ & $\Omega_2^6$\\[3pt]
          &$-d^{-2}(n+1) \left(3d^{n+1} + \alpha^2d^{n-1}(n-1) + \beta_2^2d^{n-1}(n-1) + \gamma_1^2d^{n-1}(n-1)\right)$ & $\Omega_2^7$\\[3pt]
          &$-d^{-2}(n+1) \left(2d^{n+1} + \alpha^2d^{n-1}(n-1) + \beta_2^2d^{n-1}(n-1)\right)$ & $\Omega_2^8$\\[3pt]
          &$-d^{-2}(n+1) \left(3d^{n+1} + \alpha^2d^{n-1}(n-1) + \beta_2^2d^{n-1}(n-1) + \gamma_2^2d^{n-1}(n-1)\right)$ & $\Omega_2^9$\\
          \midrule
          $\lambda_{1}=$ & $\omega$ & $\Gamma_{1}$ \\
          $\lambda_{2}=$ & $\omega$ & $\Gamma_{2}$ \\
          \midrule
          $p_2 =$ & $0$ & $\partial_1\Omega_2$ \\
          $p_2 =$ & $0$ & $\partial_2\Omega_2$ \\
          \midrule
          $p_1=$ & $-\omega$ & $\Omega_1$ \\
          \midrule
          $\vecu_1=$ & $\left[\begin{matrix}0 & 2\gamma_1^2\gamma_2^2(\beta_1\beta_2^2 + \beta_1^2\beta_2) & 2\beta_1^2\beta_2^2(\gamma_1\gamma_2^2 + \gamma_1^2\gamma_2) \end{matrix}\right]$ & $\Omega_1$ \\
          \midrule
          $\sum_{j\in\hat{S}_1}\lambda_j=$ & $2\omega$ & $\Omega_1$ \\
          \midrule
          $f_1=$ & $\beta_{1}^{2} \gamma_{2}^{2} + 4 \beta_{1} \beta_{2} \gamma_{2}^{2} + \beta_{2}^{2} \gamma_{1}^{2} + 4 \beta_{2}^{2} \gamma_{1} \gamma_{2} + 2 \beta_{2}^{2} \gamma_{2}^{2} - 2\omega$  & $\Omega_1$ \\
          \bottomrule
    \end{tabular}
    \end{adjustbox}
\end{table}

\end{appendices}

\end{document}